%% file: main.tex
\pdfoutput=1
\documentclass[onefignum,onetabnum]{siamonline250211}

\input{shared}

\ifpdf
\hypersetup{
  pdftitle={Degree-stratified randomized projection for spectral clustering},
  pdfauthor={J. Eom, Y. Hwang, E. Kim, and Y. Jang}
}
\fi

\begin{document}

\maketitle

% REQUIRED
\begin{abstract}
Many real-world networks have heavy-tailed degree distributions, in which a few high-degree nodes coexist with many low-degree nodes. Randomized range finders reduce the cost of computing the eigenspace used in spectral clustering, but their random test matrices are constructed independently of the node degrees. We propose degree-stratified randomized projection (DS-RP), which partitions the nodes into dyadic degree strata, distributes a fixed sketch dimension across the strata, and draws an independent local sketch for each stratum. We bound the angle between the sampled range and the target eigenspace and show that the gain comes from the allocation of the sketch dimension rather than from the stratification itself. The bound improves when more sketch dimensions are assigned to the strata on which the target eigenspace is concentrated. For the degree-corrected stochastic block model, we prove that the leverage scores of this eigenspace are proportional to the node degrees, so that the allocation can be determined from the degrees alone. A Rayleigh--Ritz bound that requires no smallness assumption carries the estimate over to the computed spectral embedding. Numerical experiments on synthetic benchmarks and real-world networks show that DS-RP improves the clustering accuracy of standard Gaussian and CountSketch projections at the same sketch dimension and essentially the same computational cost.
\end{abstract}

% REQUIRED

% REQUIRED

%%%%%%%%%%%%%%%%%%%%%%%%%%%%%%%%%%%%%%%%%%%%%%%%%%
\section{Introduction}\label{sec:intro}

Spectral clustering is a widely used approach for identifying community structure in graphs and networks \cite{ng2002spectral,vonluxburg2007tutorial}. It constructs a low-dimensional representation from selected eigenvectors of a graph operator and then applies a clustering algorithm, typically $k$-means, to the embedded nodes. For large-scale graphs, computing the required eigenspace often constitutes the dominant computational cost.

Randomized numerical linear algebra provides an attractive way to reduce this cost \cite{halko2011finding,mahoney2011randomized,woodruff2014sketching,martinsson2020randomized}. Randomized range finders approximate a dominant subspace by multiplying the target matrix by a low-dimensional random test matrix, orthogonalizing the sample, and solving a much smaller projected problem. Subspace iteration and block Krylov variants improve the accuracy when the spectrum decays slowly \cite{musco2015randomized,saibaba2019randomized}. These ideas have been incorporated into spectral clustering in several ways: random projection and sampling of the data \cite{sakai2009fast}, preliminary data compression \cite{yan2009fast}, a randomized power method with provable guarantees for the subsequent $k$-means step \cite{boutsidis2015spectral}, graph filtering of random signals combined with node sampling \cite{tremblay2016compressive}, and random projection and random sampling analyzed statistically under the stochastic block model and its degree-corrected variant \cite{zhang2022randomized}.

Many real-world networks exhibit degree heterogeneity, in which a large population of low-degree nodes coexists with a small number of high-degree nodes. At the modeling level, degree heterogeneity is addressed by the degree-corrected stochastic block model (DC-SBM) \cite{karrer2011stochastic} and by degree-aware spectral methods \cite{qin2013regularized,jin2015fast}. At the sketching level, in contrast, conventional randomized range finders use a single global test matrix whose construction is independent of the node degrees. Nonuniform constructions in randomized numerical linear algebra are usually guided by leverage scores, which measure how strongly each coordinate participates in a target subspace \cite{drineas2008relative,mahoney2011randomized}. Computing them, however, requires the very eigenspace that one wishes to approximate. This raises the question of whether the degree sequence, which is available at negligible cost, can guide the design of the test matrix.

In this paper, we show that it can. The idea is to partition the nodes into dyadic degree strata, to distribute a fixed total sketch dimension across the strata, and to draw an independent local sketch for each stratum. Our contributions are as follows.
\begin{itemize}[leftmargin=*]
\item We propose a new algorithm, \emph{degree-stratified randomized projection} (DS-RP), which adapts the test matrix to the degree sequence at the cost of the standard randomized range finder.
\item We derive an error bound for the sampled range and the computed embedding, which shows that the gain comes from the allocation of the sketch dimension, not from the stratification.
\item We prove that leverage scores are proportional to degrees under the DC-SBM, which justifies an allocation based on degrees alone.
\item We demonstrate on synthetic and real-world networks that DS-RP improves the clustering accuracy of standard Gaussian and CountSketch projections at the same cost.
\end{itemize}

The paper is organized as follows. \Cref{sec:preliminary} reviews spectral clustering and randomized range finding. \Cref{sec:method} introduces the proposed method, and \Cref{sec:analysis} presents its analysis. \Cref{sec:numerics} reports the numerical experiments, and \Cref{sec:conclusion} concludes the paper.

%%%%%%%%%%%%%%%%%%%%%%%%%%%%%%%%%%%%%%%%%%%%%%%%%%
\section{Preliminaries}\label{sec:preliminary}

\subsection{Spectral clustering}\label{ssec:spectral}

Let $\mathcal{G}=(V,E)$ be a connected undirected graph with node set $V=\{v_1,\ldots,v_n\}$ and adjacency matrix $A\in\mathbb{R}^{n\times n}$. The degree of node $v_i$ is $d_i=(A\bs{1})_i>0$, where $\bs{1}$ denotes the vector of all ones, and $D=\diag(d_1,\ldots,d_n)$ is the degree matrix. The normalized graph Laplacian is $\mathcal{L}=I-D^{-1/2}AD^{-1/2}$, and we define the normalized adjacency matrix
\begin{equation}
\label{eq:normalized-adjacency}
S = D^{-1/2} A D^{-1/2} = I-\mathcal{L}.
\end{equation}
The matrices $S$ and $\mathcal{L}$ share the same eigenvectors, while the ordering of their eigenvalues is reversed. We formulate spectral clustering in terms of $S$. Let $(\lambda_i,\bs{u}_i)$, $i=1,\ldots,n$, denote the eigenpairs of $S$, where the eigenvectors are orthonormal and
\[
1=\lambda_1\geq\lambda_2\geq\cdots\geq\lambda_n\geq -1 .
\]
The leading eigenvector $\bs{u}_1=D^{1/2}\bs{1}/\|D^{1/2}\bs{1}\|_2$ carries no cluster information and is excluded from the spectral embedding. The embedding is constructed from the eigenvectors associated with the $k$ largest nontrivial eigenvalues of $S$,
\[
Z_k=[\bs{u}_2,\ldots,\bs{u}_{k+1}]\in\mathbb{R}^{n\times k},
\]
where the $i$th row of $Z_k$ represents node $v_i$. Each row of $Z_k$ is normalized to unit Euclidean length, and the $k$-means algorithm \cite{lloyd1982least} is applied to the rows of the normalized embedding to obtain a cluster label for each node \cite{ng2002spectral,vonluxburg2007tutorial}.

\subsection{Randomized range finder}\label{ssec:rrf}

Randomized range finders approximate the dominant eigenspace of a matrix from its action on a low-dimensional random test matrix~\cite{halko2011finding,martinsson2020randomized}. For a symmetric matrix $H\in\mathbb{R}^{n\times n}$ and target rank $k$, let $r\geq 0$ be an oversampling parameter, set $\ell=k+r$, and let $\Omega\in\mathbb{R}^{n\times\ell}$ be a Gaussian test matrix with independent $\mathcal{N}(0,1/\ell)$ entries. Starting from $Y=H\Omega$, we perform $q\geq 0$ steps of stabilized subspace iteration,
\[
Y\leftarrow\orth(H^2Y),
\]
where $\orth(\cdot)$ denotes an orthonormal basis for the column space of its argument. The resulting matrix spans the range of $H^{2q+1}\Omega$. Larger values of $q$ improve the separation of the dominant eigenspace when the spectrum decays slowly, at the expense of additional matrix multiplications. We then compute $Q=\orth(Y)\in\mathbb{R}^{n\times\ell}$ and the eigendecomposition of the projected matrix,
\[
C=Q^\top H Q=\widetilde U\diag(\theta_1,\ldots,\theta_\ell)\widetilde U^\top,
\qquad \theta_1\geq\cdots\geq\theta_\ell,
\]
where $\theta_i$ are the Ritz values of $H$ with respect to $\range(Q)$. The Ritz vectors $U_k=Q\widetilde U_{:,1:k}\in\mathbb{R}^{n\times k}$ approximate the leading eigenvectors of $H$. For a sparse matrix $H$, the cost of this procedure is $O\bigl((2q+2)\nnz(H)\ell+(q+2)n\ell^2+\ell^3\bigr)$.

Two aspects of this procedure matter for spectral clustering. First, the iteration amplifies each eigendirection according to the absolute value of its eigenvalue, so that large negative eigenvalues of an indefinite matrix compete with the desired algebraically largest ones. Second, when applied to $S$, the trivial eigenvector $\bs{u}_1$ is the most dominant direction, so that part of the sketch is spent on a vector that is known in closed form and carries no cluster information. Since $\bs{u}_1$ is available from the degrees alone, it can be removed exactly and at negligible cost by deflation with the orthogonal projector
\begin{equation}
\label{eq:orthogonal-projector}
P_\perp=I-\bs{u}_1\bs{u}_1^\top ,
\end{equation}
that is, by applying the range finder to $P_\perp SP_\perp=S-\bs{u}_1\bs{u}_1^\top$. The computed basis $Q$ and the Ritz vectors are then exactly orthogonal to $\bs{u}_1$. The first issue is addressed by the shift introduced in \Cref{ssec:operator}. Random projection-based spectral clustering normalizes the rows of $U_k$ and applies $k$-means to the normalized embedding. Random sampling-based variants instead sparsify the graph matrix before computing a low-rank approximation~\cite{zhang2022randomized}.

\subsection{Degree heterogeneity}
\label{ssec:degree-heterogeneity}

\begin{figure}[t]
    \centering
    \includegraphics[width=\textwidth]
    {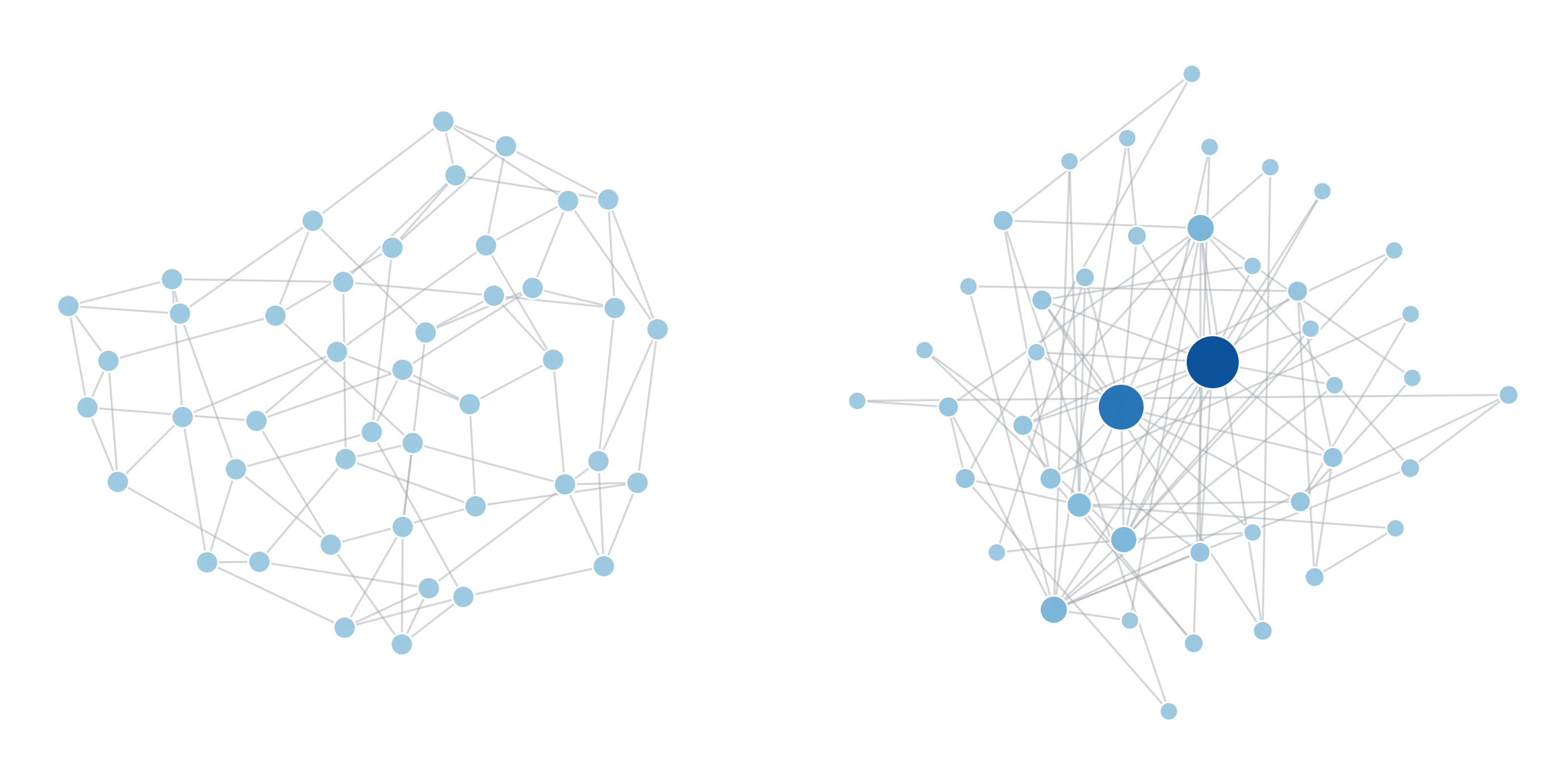}
    \caption{
        Comparison of degree-homogeneous and degree-heterogeneous graphs with comparable mean degrees. The latter contains many low-degree nodes and a small number of high-degree hubs.
    }
    \label{fig:degree-heterogeneity}
\end{figure}

Degree heterogeneity refers to the extent to which the degrees $d_1,\ldots,d_n$ vary across the nodes. As illustrated in \Cref{fig:degree-heterogeneity}, graphs with similar mean degrees can exhibit markedly different degree profiles. In particular, many low-degree nodes may coexist with a small number of high-degree hubs. A power-law degree distribution, in which the number of nodes with a given degree decreases polynomially as the degree increases, is a typical example \cite{clauset2009power}.

%%%%%%%%%%%%%%%%%%%%%%%%%%%%%%%%%%%%%%%%%%%%%%%%%%
\section{Degree-stratified randomized spectral clustering}
\label{sec:method}

The main idea of the proposed method is to replace the single global test matrix of the randomized range finder by a collection of local sketching matrices supported on degree strata. The total sketch dimension is kept fixed, but it is explicitly distributed across the degree regimes.

\subsection{Shifted and deflated graph operator}
\label{ssec:operator}

As noted in \Cref{ssec:rrf}, large negative eigenvalues of $S$ would compete with the desired ones in the subspace iteration. To avoid this, we use the shifted operator and its deflation,
\begin{equation}
\label{eq:deflated-operator}
T := \frac{1}{2}(I+S),
\qquad
T_\perp:=P_\perp T P_\perp .
\end{equation}

\begin{proposition}
\label{prop:shift-deflation}
Let $\mu_i=\frac12(1+\lambda_i)$ for $i=1,\ldots,n$. Then $T\bs{u}_i=\mu_i\bs{u}_i$ for all $i$, the eigenvalues satisfy $1=\mu_1\geq\mu_2\geq\cdots\geq\mu_n\geq 0$, and
\[
    T_\perp \bs{u}_1=0,
    \qquad
    T_\perp \bs{u}_i=\mu_i\bs{u}_i\quad\text{for } i\geq 2 .
\]
\end{proposition}

Hence $T_\perp$ is positive semidefinite, removes the trivial eigendirection, and preserves the nontrivial eigenvectors of $S$ together with the ordering of their eigenvalues. The proof is given in \cref{supp:method}.

\subsection{Degree-stratified sketch construction}
\label{ssec:ds-sketch}

We first partition the nodes into dyadic degree strata. Let $d_{\max}=\max_{1\le i\le n}d_i$ and $s_0=1+\lfloor\log_2 d_{\max}\rfloor$, and define
\[
\mathcal B_h=\left\{i\in\{1,\ldots,n\}:2^{h-1}\le d_i<2^h\right\},
\qquad h=1,\ldots,s_0 .
\]
We remove the empty sets and relabel the remaining sets as $\mathcal B_1,\ldots,\mathcal B_s$ without changing their order, so that $s$ is the number of active strata.

Given the total sketch dimension $\ell$, DS-RP assigns a local sketch dimension $\ell_j$ to each stratum, with $\sum_{j=1}^{s}\ell_j=\ell$. Each stratum receives at least $\ell_{\min}$ dimensions, which ensures that every active stratum contributes to the sampled range. Unless otherwise stated, we set $\ell_{\min}=1$. If $s\ell_{\min}>\ell$, the adjacent pair $(\mathcal B_j,\mathcal B_{j+1})$ with the smallest combined degree volume is merged repeatedly until $s\ell_{\min}\leq\ell$. For the resulting strata, let
\[
    m_j=\sum_{i\in\mathcal B_j}d_i,
    \qquad j=1,\ldots,s,
\]
denote the degree volumes. The remaining $\ell-s\ell_{\min}$ dimensions are distributed according to square-root degree-volume weights,
\begin{equation}
\label{eq:rho}
    \rho_j
    =
    (\ell-s\ell_{\min})
    \frac{\sqrt{m_j}}{\sum_{h=1}^{s}\sqrt{m_h}} .
\end{equation}
We add $\lfloor\rho_j\rfloor$ dimensions to stratum $j$ and assign the dimensions that are still left in decreasing order of the fractional parts $\rho_j-\lfloor\rho_j\rfloor$. The degree volume $m_j$ serves as a computable proxy for the leverage that stratum $\mathcal B_j$ carries for the target eigenspace, and the square root tempers the concentration of the budget in a few high-volume strata. Both points are discussed in \Cref{sec:analysis}.

For each stratum, let $\Omega_j\in\mathbb R^{|\mathcal B_j|\times \ell_j}$ be an independent local sketching matrix. The stratum-wise samples are concatenated,
\[
    Y_j=(T_\perp)_{:,\mathcal B_j}\Omega_j,
    \qquad
    Y=[Y_1,\ldots,Y_s]\in\mathbb R^{n\times \ell}.
\]
We consider two families of local sketches, which use the same strata and the same allocation. In the Gaussian variant, the entries of $\Omega_j$ are independent $\mathcal N(0,1/\ell_j)$ random variables. In the CountSketch variant \cite{clarkson2017low}, each row of $\Omega_j$ has a single nonzero entry, a random sign placed in a column chosen uniformly at random.

The local sketching matrices are drawn only once. After $Y$ has been formed, the algorithm proceeds as in \Cref{ssec:rrf} with $H=T_\perp$. It applies $q$ steps of subspace iteration, $Y \leftarrow \orth(T_\perp^2Y)$, computes $Q=\orth(Y)$, and extracts the $k$ leading Ritz vectors $U_k$ of $C=Q^\top T_\perp Q$, which provide the approximate nontrivial eigenvectors used for the spectral embedding. The complete procedure is summarized in \Cref{alg:degree-stratified-rp}.

\begin{algorithm}[!htp]
\caption{Degree-stratified randomized projection (DS-RP)}
\label{alg:degree-stratified-rp}
% \small
\begin{algorithmic}[1]
\Require Symmetric adjacency matrix
$A\in\mathbb{R}^{n\times n}$ with positive degrees,
target dimension $k$, oversampling parameter $r$,
number of subspace iterations $q$, minimum allocation $\ell_{\min}$,
and a sketching distribution
\Ensure Approximate nontrivial eigenvectors $U_k$ and eigenvalues
$\Lambda_k$ of $S$

\State $\ell\gets k+r$, \quad
$\bs{d}\gets A\bs{1}$, \quad
$D\gets\operatorname{diag}(\bs{d})$, \quad
$S\gets D^{-1/2}AD^{-1/2}$
\State $T\gets\frac12(I+S)$
\State $\bs{u}_1\gets
D^{1/2}\bs{1}/\|D^{1/2}\bs{1}\|_2$
\State $P_\perp\gets I-\bs{u}_1\bs{u}_1^\top$, \quad
$T_\perp\gets P_\perp T P_\perp$

\State $d_{\max}\gets\max_{1\leq i\leq n}d_i$, \quad
$s_0\gets 1+\lfloor\log_2 d_{\max}\rfloor$
\State $\mathcal{B}_h\gets\{i\in\{1,\ldots,n\}:2^{h-1}\le d_i<2^h\}$
for $h=1,\ldots,s_0$
\State Relabel the nonempty sets as $\mathcal{B}_1,\ldots,\mathcal{B}_s$

\State Merge the adjacent pair
$(\mathcal{B}_j,\mathcal{B}_{j+1})$ minimizing
\[
    \sum_{i\in\mathcal{B}_j\cup\mathcal{B}_{j+1}} d_i
\]
until $s\ell_{\min}\leq\ell$, updating $s$ and relabeling after each merge

\State $m_j\gets\displaystyle\sum_{i\in\mathcal{B}_j}d_i$ for $j=1,\ldots,s$
\State $\ell_j\gets\ell_{\min}$ for $j=1,\ldots,s$
\State $\rho_j\gets(\ell-s\ell_{\min})\sqrt{m_j}\big/\sum_{h=1}^{s}\sqrt{m_h}$, \quad
$\ell_j\gets\ell_j+\lfloor\rho_j\rfloor$
for $j=1,\ldots,s$
\State Assign the remaining
$\ell-\sum_{j=1}^{s}\ell_j$ dimensions in decreasing order of
$\rho_j-\lfloor\rho_j\rfloor$

\For{$j=1,\ldots,s$}
    \State Draw a sketching matrix
    $\Omega_j\in\mathbb{R}^{|\mathcal{B}_j|\times\ell_j}$
    from the sketching distribution
    \State $Y_j\gets
    (T_\perp)_{:,\mathcal{B}_j}\Omega_j$
\EndFor

\State $Y\gets[Y_1,\ldots,Y_s]$

\State Repeat $q$ times: $Y\gets\operatorname{orth}(T_\perp^2Y)$

\State $Q\gets\operatorname{orth}(Y)$
\State $C\gets Q^\top T_\perp Q$
\State Compute the eigendecomposition
\[
C
=
\widetilde U
\operatorname{diag}(\theta_1,\ldots,\theta_\ell)
\widetilde U^\top,
\qquad
\theta_1\geq\cdots\geq\theta_\ell
\]
\State $U_k\gets Q\widetilde U_{:,1:k}$
\State $\Lambda_k\gets
\operatorname{diag}(2\theta_1-1,\ldots,2\theta_k-1)$
\State \Return $U_k,\Lambda_k$

\end{algorithmic}
\end{algorithm}

DS-RP has the same leading-order cost as a standard randomized range finder with the same sketch dimension. The stratification and the allocation require only the degree vector. Neither $P_\perp$ nor $T_\perp$ is formed explicitly. Applying $T_\perp$ to an $n\times \ell$ matrix requires one multiplication with $S$ and two projections onto the orthogonal complement of $\bs{u}_1$, which takes $O(\nnz(A)\ell+n\ell)$ operations. Since the algorithm performs $2q+2$ such multiplications, the total cost of constructing the spectral embedding is $O\bigl((2q+2)\nnz(A)\ell+(q+2)n\ell^2+\ell^3\bigr)$.

%%%%%%%%%%%%%%%%%%%%%%%%%%%%%%%%%%%%%%%%%%%%%%%%%%
\section{Analysis of degree-stratified randomized projection}
\label{sec:analysis}

We analyze the sketch constructed in \Cref{sec:method} for Gaussian local sketches. All proofs are given in \cref{supp:analysis}.

\subsection{Range capture with a stratified test matrix}
\label{ssec:structural}

By \Cref{prop:shift-deflation}, the eigenvalues of $T_\perp$ are $\mu_2\geq\cdots\geq\mu_n\geq 0$, with eigenvectors $\bs{u}_2,\ldots,\bs{u}_n$, together with the eigenvalue $0$ associated with $\bs{u}_1$. \Cref{alg:degree-stratified-rp} extracts $k$ Ritz vectors, but the eigenvectors that carry the community structure may span a smaller space. For example, for the DC-SBM analyzed in \Cref{prop:dcsbm_leverage}, they are $\bs{u}_2,\ldots,\bs{u}_k$. We therefore fix a \emph{signal dimension} $k_\star\leq k$ at which the spectrum of $T_\perp$ has a gap, $\gamma_\star:=\mu_{k_\star+1}-\mu_{k_\star+2}>0$, and write $Z_\star=[\bs{u}_2,\ldots,\bs{u}_{k_\star+1}]$ for the target eigenvectors.

Let $n_j=|\mathcal B_j|$. If the nodes are ordered by strata, the sample matrix of \Cref{alg:degree-stratified-rp} is
\begin{equation}
\label{eq:structured-test-matrix}
    Y=T_\perp\Omega,
    \qquad
    \Omega:=
    \begin{bmatrix}
    \Omega_1 & & \\
    & \ddots & \\
    & & \Omega_s
    \end{bmatrix}
    \in\mathbb R^{n\times\ell}.
\end{equation}
Hence DS-RP is a randomized range finder with a block-diagonal test matrix, and the standard randomized range finder corresponds to a single stratum, $s=1$. We measure how well the sampled range captures the target eigenspace by
\[
    \dist(Z_\star,Q):=\norm{(I-QQ^\top)Z_\star}_2,
\]
the sine of the largest principal angle between $\range(Z_\star)$ and $\range(Q)$.

The accuracy of the sampled range depends on how the target eigenspace is distributed over the strata. Let $\bs{z}_i\in\mathbb R^{k_\star}$ denote the $i$th row of $Z_\star$, written as a column vector. The \emph{leverage score} of node $v_i$ with respect to the target eigenspace is $\zeta_i=\norm{\bs{z}_i}_2^2$ \cite{drineas2008relative,mahoney2011randomized}. For each stratum, we define
\[
    W_j=\sum_{i\in\mathcal B_j}\bs{z}_i\bs{z}_i^\top,
    \qquad
    \tau_j=\trace(W_j)=\sum_{i\in\mathcal B_j}\zeta_i ,
\]
and call $W_j$ the \emph{stratum leverage matrix} and $\tau_j$ the \emph{stratum leverage} of $\mathcal B_j$. Since the strata partition the node set, $\sum_{j}W_j=Z_\star^\top Z_\star=I$ and $\sum_{j}\tau_j=k_\star$. The allocation enters the analysis through the quantity
\begin{equation}
\label{eq:design-functional}
    \Phi:=\lambda_{\min}\Biggl(\sum_{j=1}^{s}\frac{\ell_j}{n_j}\,W_j\Biggr),
\end{equation}
where $\lambda_{\min}(\cdot)$ denotes the smallest eigenvalue of a symmetric matrix.

\begin{theorem}
\label{thm:ds_gaussian}
Let $Q$ be computed by \Cref{alg:degree-stratified-rp} with Gaussian local sketches. Let $t>0$, and suppose that $\ell_j\leq n_j$ and $\sqrt{\ell_j}\geq 2(\sqrt{k_\star}+t)$ for $j=1,\ldots,s$. Then, with probability at least $1-2s\,e^{-t^2/2}$,
\[
    \dist(Z_\star,Q)
    \leq
    5\left(\frac{\mu_{k_\star+2}}{\mu_{k_\star+1}}\right)^{2q+1}
    \frac{1}{\sqrt{\Phi}} .
\]
\end{theorem}

The factor $(\mu_{k_\star+2}/\mu_{k_\star+1})^{2q+1}$ is the usual contraction of subspace iteration and does not depend on the test matrix, which enters only through $\Phi$. For a single stratum, $W_1=I$ and $\Phi=\ell/n$, which is the familiar dependence of the Gaussian range finder on the sketch dimension \cite{halko2011finding}. Since $\sum_jW_j=I$, the proportional allocation $\ell_j=\ell n_j/n$ gives $\Phi=\ell/n$ as well, so that stratification by itself brings no gain. The bound improves on that of the standard sketch whenever the allocation satisfies $\Phi>\ell/n$, that is, whenever sketch dimensions are shifted toward strata that carry more leverage per node while every direction of $\range(Z_\star)$ remains covered. Choosing the allocation that maximizes $\Phi$ is an E-optimal design problem \cite{pukelsheim2006optimal}.

The choice of the error measure matters for this conclusion. It is tempting to assess a sketch instead by the accuracy of $YY^\top$, formed from the initial sample matrix $Y=T_\perp\Omega$, as an estimator of $T_\perp^2$. For Gaussian local sketches this estimator is unbiased, but for every stratification and every allocation its mean-square error is at least as large as that of the standard Gaussian sketch (\cref{prop:gram_mse}). This criterion measures the approximation error over the entire operator, including spectral components outside the target eigenspace, whereas spectral clustering depends only on $\range(Z_\star)$.

\subsection{Leverage and degree}
\label{ssec:allocation}

The stratum leverage matrices $W_j$ depend on the unknown eigenvectors. We now show that, under the DC-SBM \cite{karrer2011stochastic}, leverage is determined by degree. Consider a population graph with communities $\mathcal C_1,\ldots,\mathcal C_k$ and weighted adjacency matrix $\mathcal A$ with entries
\begin{equation}
\label{eq:dcsbm}
    \mathcal A_{ii'}=\psi_i\psi_{i'}B_{ab}
    \qquad\text{for } v_i\in\mathcal C_a \text{ and } v_{i'}\in\mathcal C_b,
\end{equation}
where $\psi_1,\ldots,\psi_n$ are positive degree-correction factors and $B\in\mathbb R^{k\times k}$ is symmetric with nonnegative entries. For a set of nodes $\mathcal C$, let $\vol(\mathcal C)=\sum_{v_i\in\mathcal C}d_i$, and let $\vol(\mathcal G)=\sum_{i=1}^{n}d_i$.

\begin{proposition}
\label{prop:dcsbm_leverage}
Let $\mathcal A$ be given by \eqref{eq:dcsbm} with $B$ nonsingular, let the degrees $d_i=(\mathcal A\bs{1})_i$ be positive, and let $S=D^{-1/2}\mathcal A D^{-1/2}$. Then $S$ has rank $k$, and the eigenvectors associated with its nonzero eigenvalues other than $\lambda_1=1$ span a subspace of dimension $k-1$ whose leverage scores are
\[
    \zeta_i=d_i\left(\frac{1}{\vol(\mathcal C(v_i))}-\frac{1}{\vol(\mathcal G)}\right),
    \qquad i=1,\ldots,n,
\]
where $\mathcal C(v_i)$ denotes the community that contains $v_i$. In particular, if all communities have the same volume, then $\zeta_i=(k-1)\,d_i/\vol(\mathcal G)$.
\end{proposition}

If $B$ is positive definite, these eigenvalues are positive, so that the subspace in \Cref{prop:dcsbm_leverage} is the target eigenspace $\range(Z_\star)$ with $k_\star=k-1$. Summing the leverage scores over a stratum then gives
\begin{equation}
\label{eq:tau-volume}
    \tau_j=(k-1)\,\frac{m_j}{\vol(\mathcal G)} ,
    \qquad j=1,\ldots,s,
\end{equation}
exactly for the population model with communities of equal volume, and approximately for graphs sampled from it. The stratum leverage is thus proportional to the degree volume $m_j$ used in \Cref{alg:degree-stratified-rp}. If, moreover, the leverage of each stratum is spread evenly over the directions of the target eigenspace, $W_j\approx(\tau_j/k_\star)I$, then \eqref{eq:design-functional} and \eqref{eq:tau-volume} give
\begin{equation}
\label{eq:gain}
    \Phi\approx\frac{1}{\vol(\mathcal G)}\sum_{j=1}^{s}\ell_j\,\frac{m_j}{n_j} ,
\end{equation}
so that each sketch dimension contributes the mean degree $m_j/n_j$ of the stratum to which it is assigned. The right-hand side equals $\ell/n$ for the proportional allocation and increases when sketch dimensions are moved to strata of above-average degree, which in a heavy-tailed graph contain few nodes but most of the volume.

This reallocation has limits. A stratum cannot make use of more sketch dimensions than it has nodes, which is the condition $\ell_j\leq n_j$ in \Cref{thm:ds_gaussian}, and one with fewer than $k_\star$ nodes is unable to cover all directions of $\range(Z_\star)$, which the smallest eigenvalue in \eqref{eq:design-functional} penalizes. The rule \eqref{eq:rho} is therefore a tempered version of volume-proportional allocation. It is increasing in $m_j$, as suggested by \eqref{eq:tau-volume}, while the square root and the minimum allocation $\ell_{\min}$ keep every stratum represented. The rule is a heuristic rather than an optimal design, since the allocation that maximizes $\Phi$ depends on the unknown matrices $W_j$.

\subsection{Ritz vectors and spectral embedding}
\label{ssec:eigenspace_embedding}

We finally relate the sampled range to the output of the Rayleigh--Ritz step. Let $U_k=Q\widetilde U_{:,1:k}$ be the matrix of the $k$ leading Ritz vectors returned by \Cref{alg:degree-stratified-rp}, where $k_\star\leq k\leq\ell$. The quantity $\norm{(I-U_kU_k^\top)Z_\star}_2$ is the sine of the largest principal angle between $\range(Z_\star)$ and $\range(U_k)$.

\begin{theorem}
\label{thm:rayleigh_ritz}
Let $Q\in\mathbb R^{n\times\ell}$ have orthonormal columns, and let $U_k$ be the matrix of the $k$ leading Ritz vectors of $T_\perp$ with respect to $\range(Q)$, where $k_\star\leq k\leq\ell$. Then
\[
    \norm{(I-U_kU_k^\top)Z_\star}_2
    \leq
    \sqrt{1+\frac{\mu_2^2}{4\gamma_\star^2}}\;\dist(Z_\star,Q).
\]
\end{theorem}

\Cref{thm:rayleigh_ritz} extends the classical bound for a single Ritz vector \cite[Chapter~4]{saad2011numerical} to invariant subspaces. Perturbation bounds of Davis--Kahan type \cite{davis1970rotation,stewart1990matrix} would instead compare $T_\perp$ with its compression $QQ^\top T_\perp QQ^\top$ and require the difference to be small relative to the gap $\gamma_\star$. This cannot be expected here. The compression has rank at most $\ell$, so that the norm of the difference is at least $\mu_{\ell+2}$, and for the shifted operator this bulk eigenvalue is typically much larger than $\gamma_\star$. \Cref{thm:rayleigh_ritz} requires no such assumption. The theorem also covers $k>k_\star$, where more Ritz vectors are extracted than the signal dimension. The informative eigenspace is then nearly contained in the span of the computed embedding.

\begin{corollary}
\label{cor:main}
Suppose that the assumptions of \Cref{thm:ds_gaussian} hold. Then, with probability at least $1-2s\,e^{-t^2/2}$,
\[
    \norm{(I-U_kU_k^\top)Z_\star}_2
    \leq
    5\,\sqrt{1+\frac{\mu_2^2}{4\gamma_\star^2}}\,
    \left(\frac{\mu_{k_\star+2}}{\mu_{k_\star+1}}\right)^{2q+1}
    \frac{1}{\sqrt{\Phi}} ,
\]
and there exists a matrix $R\in\mathbb R^{k\times k_\star}$ with orthonormal columns such that
\[
    \norm{U_kR-Z_\star}_F\leq\sqrt{2k_\star}\,\norm{(I-U_kU_k^\top)Z_\star}_2 .
\]
\end{corollary}

In summary, the accuracy of DS-RP is governed by three factors: the spectral gap $\gamma_\star$ of the graph operator, the contraction of subspace iteration, and the quality $\Phi$ of the test matrix. Degree stratification acts only on the last factor, by matching the support of the sketch columns to the strata that carry the leverage of the target eigenspace.

\subsection{Validation on a DC-SBM graph}
\label{ssec:illustration}

\Cref{thm:ds_gaussian} identifies $\Phi$ as the quantity that governs the error when the local sketch dimensions are sufficiently large. We now examine whether $\Phi$ remains predictive for the small sketch dimensions used in practice, on a graph sampled from the DC-SBM with $5000$ nodes and $k=4$ communities of equal size. The degree-correction factors $\psi_i$ are drawn from a log-normal distribution with shape parameter $1.1$ and clipped to $[0.05,8]$, which produces a heavy-tailed degree distribution. The target average degree is $24$, the within-to-between edge-probability ratio is $8$, and the maximum edge probability is $0.95$. The largest connected component has $n=4929$ nodes, $\num{59992}$ edges, and $s=8$ active strata. We take $k_\star=k-1=3$, $\ell=64$, $q=3$, and $\ell_{\min}=1$, and we use Gaussian local sketches.

\begin{figure}[t]
    \centering
    \includegraphics[width=0.92\textwidth]{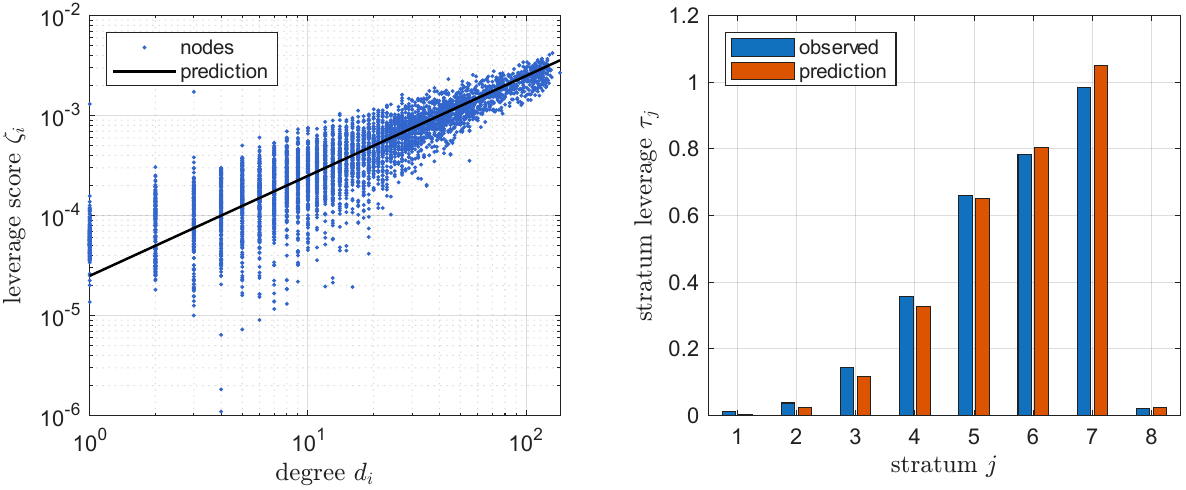}
    \caption{Leverage scores versus degrees (left) and stratum leverages (right) for the DC-SBM graph, with the predictions of \Cref{prop:dcsbm_leverage} and \eqref{eq:tau-volume}.}
    \label{fig:leverage-degree}
\end{figure}

\begin{table}[t]
\centering
\footnotesize
\caption{Average subspace errors (standard deviation) over $30$ runs for the DC-SBM graph.}
\label{tab:illustration}
\begin{tabular}{lcccc}
\toprule
Test matrix & $\Phi n/\ell$ & Predicted tangent & Observed tangent & Ritz vector error \\
\midrule
Gaussian & 1.00 & -- & 0.836 (0.060) & 0.643 (0.028) \\
DS-Gaussian, $\ell_j\propto n_j$ & 1.03 & 0.824 & 0.838 (0.074) & 0.643 (0.034) \\
DS-Gaussian, uniform & 1.39 & 0.709 & 0.646 (0.035) & 0.544 (0.021) \\
DS-Gaussian, $\ell_j\propto\sqrt{m_j}$ & 1.69 & 0.644 & 0.648 (0.052) & 0.544 (0.031) \\
DS-Gaussian, $\ell_j\propto m_j$ & 1.87 & 0.611 & 0.588 (0.057) & 0.506 (0.036) \\
\bottomrule
\end{tabular}
\end{table}

\Cref{fig:leverage-degree} shows that the leverage scores are nearly proportional to the degrees, with a correlation coefficient of $0.94$, and that the stratum leverages follow the prediction \eqref{eq:tau-volume} closely. \Cref{tab:illustration} compares the standard Gaussian sketch with DS-Gaussian sketches under four allocations: $\ell_j$ proportional to $n_j$, uniform, proportional to $\sqrt{m_j}$ as in \eqref{eq:rho}, and proportional to $m_j$. For each sketch, we report $\Phi$ relative to its value $\ell/n$ for the standard sketch, the tangent of the largest principal angle between $\range(Z_\star)$ and $\range(Q)$, to which the bound of \Cref{thm:ds_gaussian} applies as well, and the error $\norm{(I-U_kU_k^\top)Z_\star}_2$ of the Ritz vectors. The predicted tangent is the tangent observed for the standard sketch multiplied by $(\Phi n/\ell)^{-1/2}$.

The predictions agree with the observations to within about ten percent, so that $\Phi$ is predictive in this regime as well. Stratification with proportional allocation behaves like the standard sketch, as predicted, whereas the allocations with larger $\Phi$ reduce the tangent by roughly $20$ to $30$ percent. On this graph, the uniform allocation and the square-root rule perform similarly, and the volume-proportional allocation is slightly more accurate. When $q$ is varied from $0$ to $6$, the ratio between the tangents of the standard sketch and of the square-root rule stays between $1.28$ and $1.40$, in agreement with \Cref{thm:ds_gaussian}, where the factor that depends on the test matrix is independent of $q$. Finally, the error of the Ritz vectors exceeds $\dist(Z_\star,Q)$ by less than one percent in all cases, so that the Rayleigh--Ritz step loses far less than the worst-case factor in \Cref{thm:rayleigh_ritz}. The effect of the improved range capture on the clustering accuracy is examined in \Cref{sec:numerics}.

%%%%%%%%%%%%%%%%%%%%%%%%%%%%%%%%%%%%%%%%%%%%%%%%%%
\section{Numerical experiments}\label{sec:numerics}
While \Cref{ssec:illustration} examined the quantities that appear in the analysis, this section evaluates the clustering accuracy and the computational cost of the proposed DS-RP method. We compare the Gaussian and CountSketch variants of DS-RP, denoted by DS-Gaussian and DS-CountSketch (DS-CS), with standard spectral clustering and the conventional Gaussian and CountSketch (CS) methods. We refer to DS-Gaussian and DS-CS as the DS methods. All five methods are applied to the operator $T_\perp$ defined in \eqref{eq:deflated-operator}. Standard spectral clustering computes the leading eigenvectors of $T_\perp$ directly by an eigensolver and applies $k$-means to the row-normalized embedding \cite{ng2002spectral}, whereas the other methods approximate this embedding by random projection and then apply $k$-means in the same manner.
To evaluate the clustering accuracy of each method, we compute the macro-averaged F1 score (Macro-F1) after the optimal one-to-one matching between the predicted and ground-truth cluster labels \cite{weber2016comparison}, the normalized mutual information (NMI) \cite{strehl2002cluster}, and the adjusted Rand index (ARI) \cite{hubert1985comparing}. We consider synthetic benchmark graphs, graphs based on real-world data, and two large-scale graphs on which we also compare the time required to compute the embedding.

\subsection{Synthetic benchmark}\label{ssec:synthetic}
For the synthetic benchmark graphs and real-world graphs, we set $q=3$ and $r=60$. Each experiment is restricted to the largest connected component of the graph, and the target number of clusters $k$ is set to the number of distinct ground-truth labels. We perform $30$ independent runs of every method. Within a run, all methods share the same $k$-means initialization seed, while the randomized sketch of each method is generated with an independent seed. The results are averaged over the $30$ runs. For DS-RP, we set $\ell_{\min}=1$.

To evaluate the clustering performance for different degree distributions and community structures, we consider three synthetic benchmark graphs: the DC-SBM \cite{karrer2011stochastic}, the artificial benchmark for community detection (ABCD) \cite{kaminski2021artificial}, and the Lancichinetti--Fortunato--Radicchi (LFR) benchmark \cite{lancichinetti2008benchmark}.
The synthetic DC-SBM graph follows the DC-SBM with $5000$ nodes and $k=4$ latent communities. Each node $v_i$ has a degree-correction factor $\psi_i$ drawn from a log-normal distribution with shape parameter $1.1$ and clipped to $[0.05,8]$, which produces a heavy-tailed degree distribution. We set the target average degree to $24$, the within-to-between edge-probability ratio to $8$, and the maximum edge probability to $0.95$. After we restrict the graph to its largest connected component, the resulting graph has $4924$ nodes, $\num{60189}$ edges, an average degree of $24.45$, and $k=4$ ground-truth clusters, with a degree Gini coefficient of $0.5466$.
The ABCD benchmark graph follows an ABCD-style planted-community model with $5000$ nodes and $k=8$ latent communities. The target node degrees are drawn from a truncated power-law distribution with exponent $2.5$ and bounded between $10$ and $180$. The community sizes are $1852$, $874$, $583$, $448$, $371$, $322$, $287$, and $263$. The mixing parameter is $0.22$. The resulting graph is connected and has $5000$ nodes, $\num{54546}$ edges, an average degree of $21.82$, and $k=8$ ground-truth clusters, with a degree Gini coefficient of $0.3768$.
The LFR benchmark graph is generated with $5000$ nodes. The degree-distribution exponent is $2.5$, the community-size-distribution exponent is $1.5$, the mixing parameter is again $0.22$, and the average-degree parameter is $22$. The community sizes are bounded between $300$ and $1200$. This graph is also connected, with $5000$ nodes, $\num{69703}$ edges, an average degree of $27.88$, $k=9$ ground-truth clusters, and a degree Gini coefficient of $0.3599$.

\Cref{tab:accuracy-synthetic} lists the clustering accuracy of the four randomized methods on the three synthetic benchmark graphs. The Macro-F1, NMI, and ARI values obtained by standard spectral clustering, which serve as reference values, are $0.9573$, $0.8408$, and $0.8896$ for the synthetic DC-SBM graph, $0.7940$, $0.8994$, and $0.7404$ for the ABCD benchmark graph, and $0.9998$, $0.9993$, and $0.9995$ for the LFR benchmark graph, respectively. The highest value among the four randomized methods is shown in bold, and the value in parentheses gives its percentage relative to the reference value.
The DS-Gaussian and DS-CS methods give higher clustering accuracy than the conventional Gaussian and CS methods for all three graphs. In particular, the Macro-F1 score of DS-Gaussian for the ABCD benchmark graph reaches $102.29\%$ of the reference value, and for the LFR benchmark graph all three accuracy measures of DS-Gaussian exceed $99.5\%$ of their reference values. Hence, with the same sketch dimension and the same number of subspace iterations, degree stratification improves on the conventional randomized methods and recovers most of the accuracy of standard spectral clustering.
\begin{table}[htbp]
\centering
\footnotesize
\caption{Average clustering accuracy over $30$ runs on the synthetic benchmark graphs.}
\label{tab:accuracy-synthetic}
\begin{tabular}{cc cccc}
\toprule
Graph & Metric & Gaussian & CS & DS-Gaussian & DS-CS \\
\midrule
\multirow{3}{*}{Synthetic DC-SBM}
& Macro-F1  & 0.7948 & 0.7939 & {\bf 0.9013} (94.15\%) & {\bf 0.9013} (94.15\%) \\
& NMI & 0.4736 & 0.4726 & 0.6919 & {\bf 0.6922} (82.33\%) \\
& ARI & 0.5280 & 0.5267 & {\bf 0.7545} (84.81\%) & 0.7544 \\
\midrule
\multirow{3}{*}{ABCD benchmark}
& Macro-F1  & 0.7408 & 0.7402 & {\bf 0.8122} (102.29\%) & 0.7765 \\
& NMI & 0.7524 & 0.7514 & {\bf 0.8566} (95.24\%) & 0.8473 \\
& ARI & 0.6069 & 0.6061 & {\bf 0.7334} (99.05\%) & 0.6948 \\
\midrule
\multirow{3}{*}{LFR benchmark}
& Macro-F1  & 0.9577 & 0.9566 & {\bf 0.9983} (99.85\%) & 0.9979 \\
& NMI & 0.8874 & 0.8853 & {\bf 0.9948} (99.55\%) & 0.9935 \\
& ARI & 0.8889 & 0.8849 & {\bf 0.9966} (99.71\%) & 0.9956 \\
\bottomrule
\end{tabular}
\end{table}

\begin{figure}[htbp]
\centering
\begin{minipage}{0.32\linewidth}
\centering
\includegraphics[clip=true,height=1.7in]{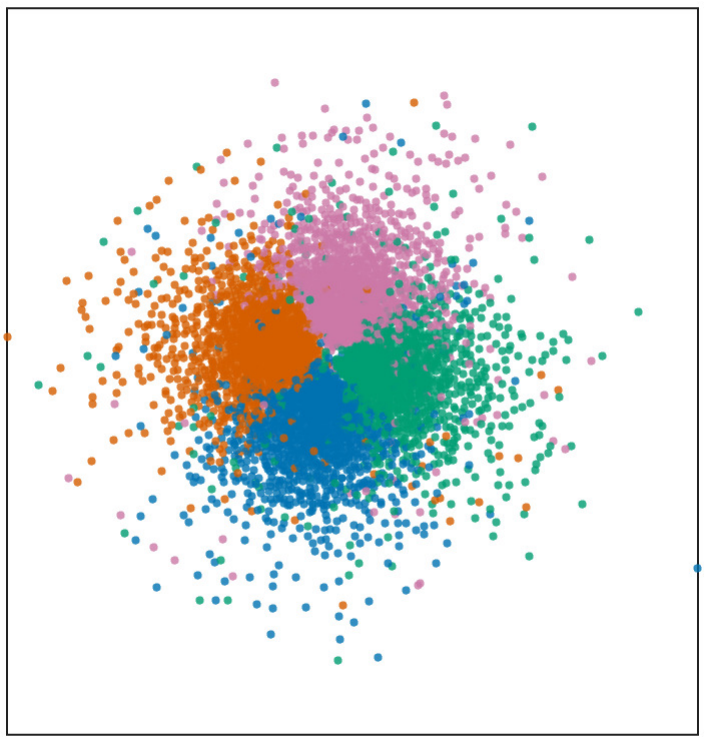}\\
Ground truth \vspace{0.3cm}
\end{minipage}
\begin{minipage}{0.32\linewidth}
\centering
\includegraphics[clip=true,height=1.7in]{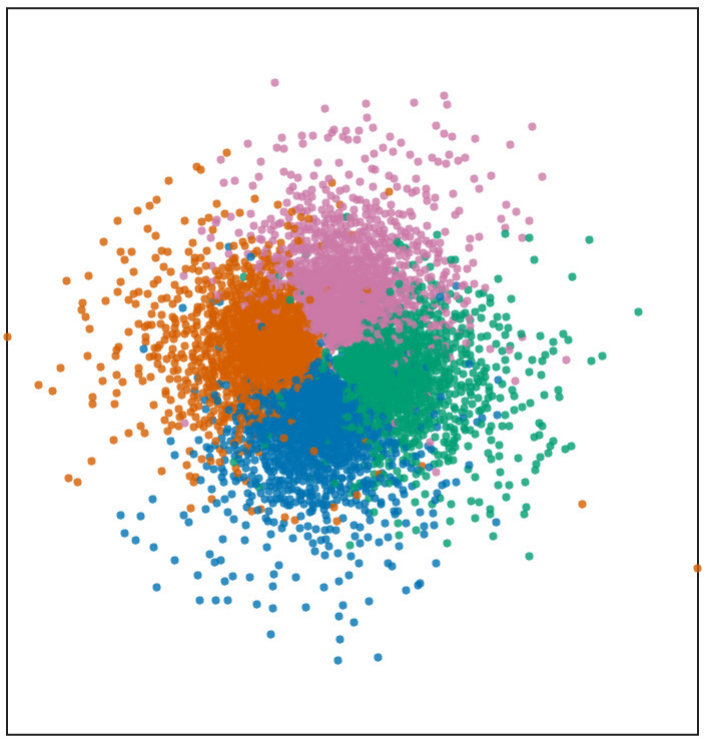}\\
Standard\vspace{0.3cm}
\end{minipage}
\begin{minipage}{0.32\linewidth}
\centering
\includegraphics[clip=true,height=1.7in]{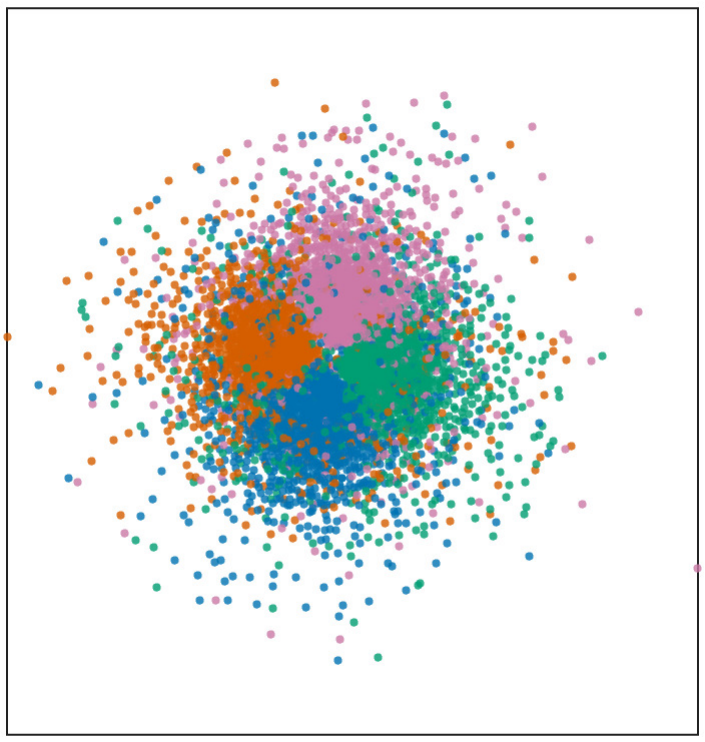}\\
Gaussian\vspace{0.3cm}
\end{minipage}\\
\begin{minipage}{0.32\linewidth}
\centering
\includegraphics[clip=true,height=1.7in]{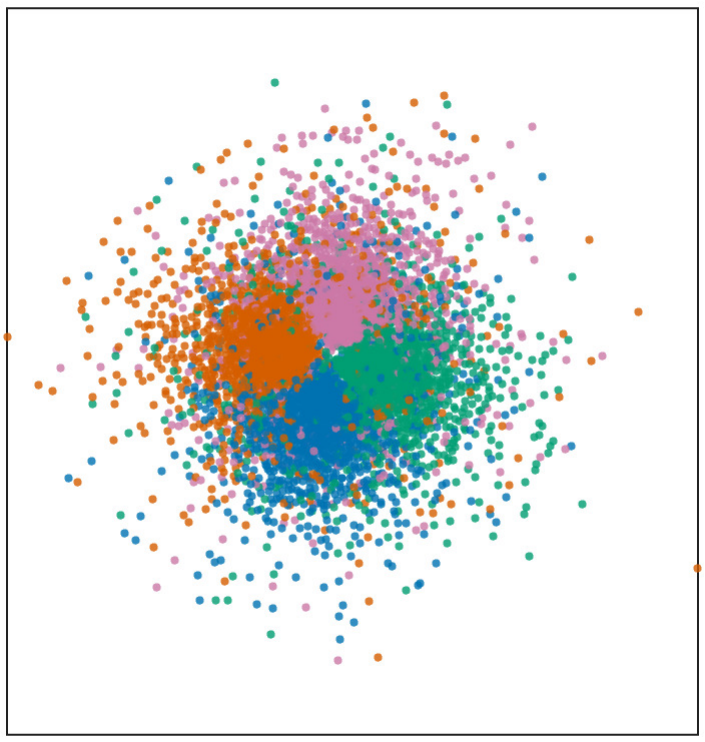}\\
CS
\end{minipage}
\begin{minipage}{0.32\linewidth}
\centering
\includegraphics[clip=true,height=1.7in]{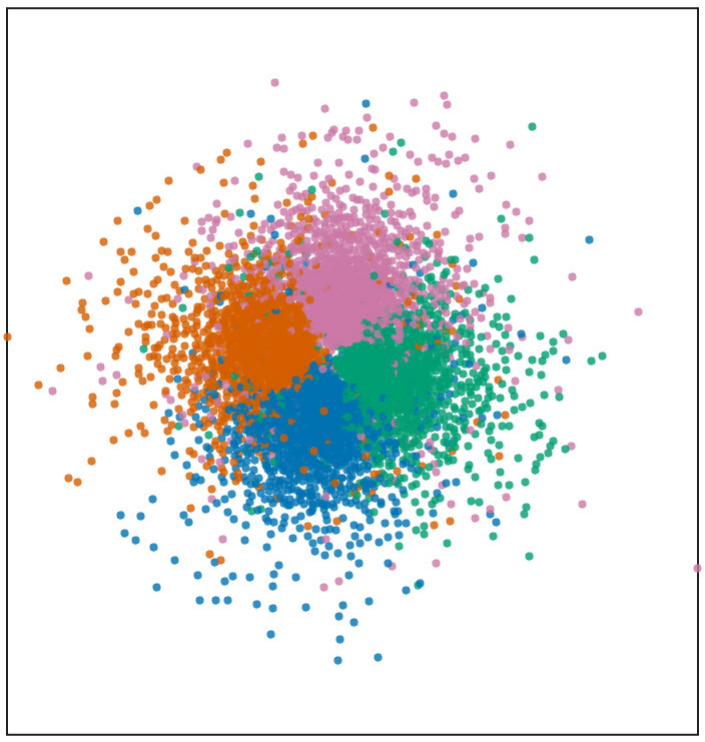}\\
DS-Gaussian
\end{minipage}
\begin{minipage}{0.32\linewidth}
\centering
\includegraphics[clip=true,height=1.7in]{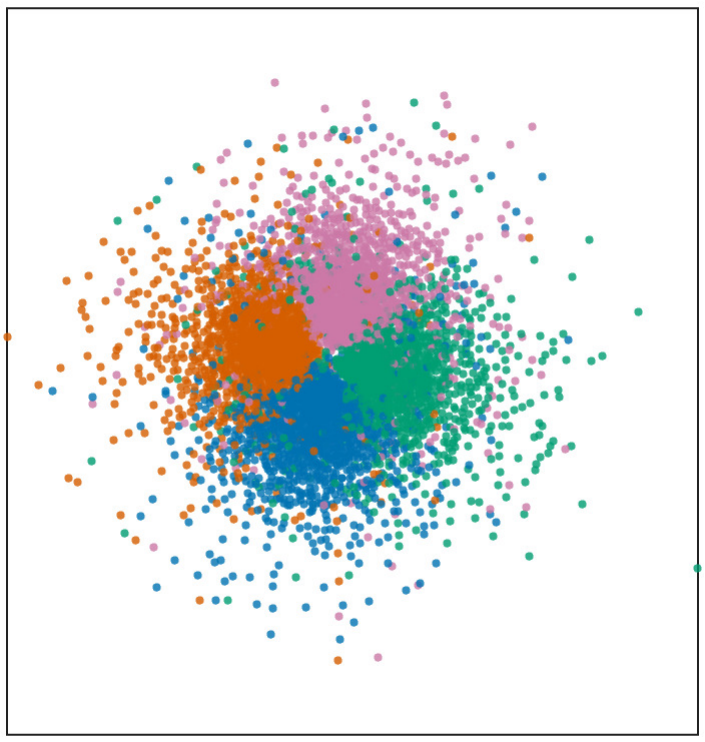}\\
DS-CS
\end{minipage}
\caption{Clustering results for the synthetic DC-SBM graph with $k=4$ clusters.}
\label{fig:synthetic-dcsbm}
\end{figure}
\begin{figure}[htbp]
\centering
\begin{minipage}{0.32\linewidth}
\centering
\includegraphics[clip=true,height=1.7in]{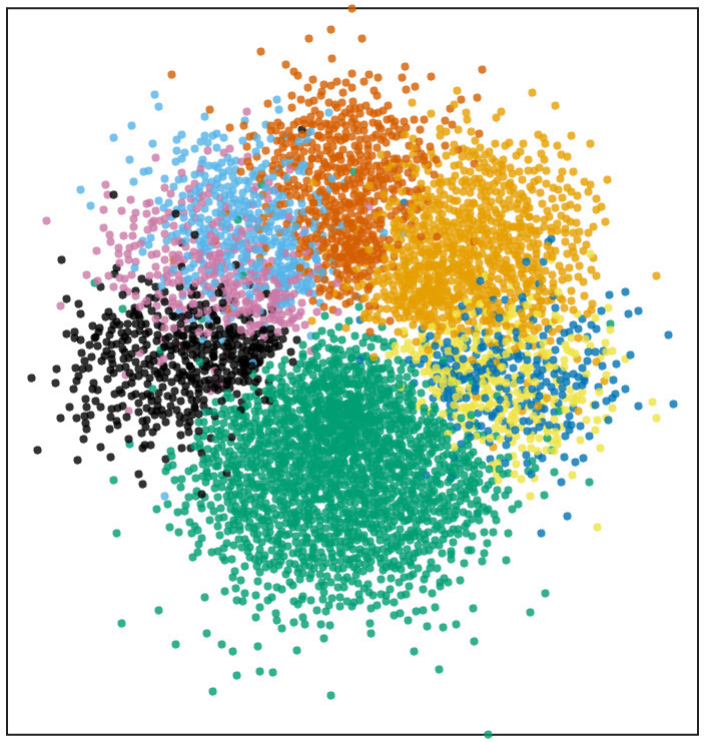}\\
Ground truth \vspace{0.3cm}
\end{minipage}
\begin{minipage}{0.32\linewidth}
\centering
\includegraphics[clip=true,height=1.7in]{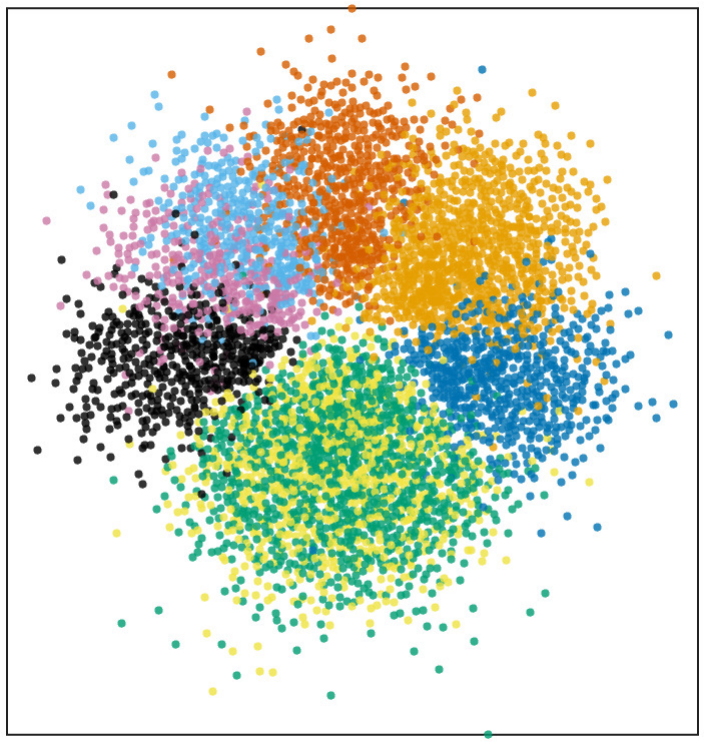}\\
Standard\vspace{0.3cm}
\end{minipage}
\begin{minipage}{0.32\linewidth}
\centering
\includegraphics[clip=true,height=1.7in]{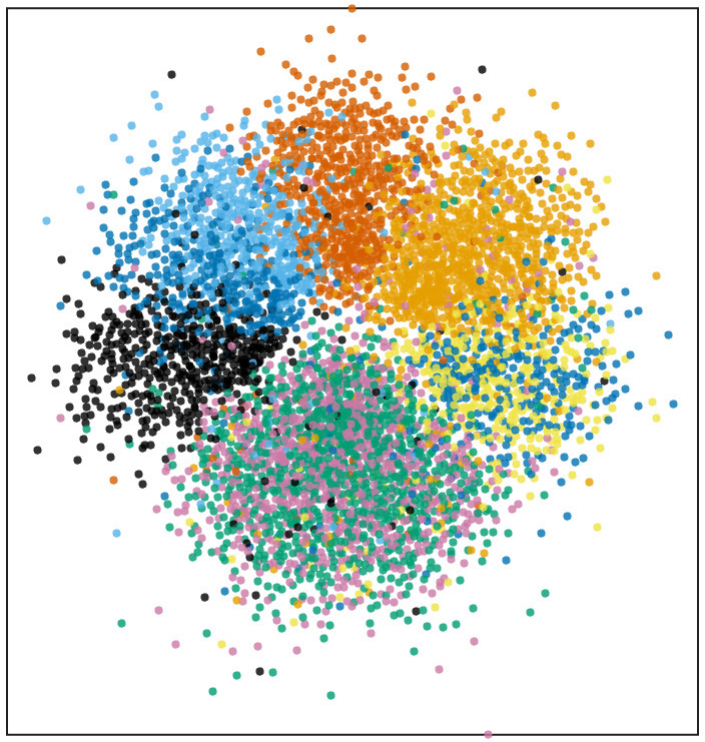}\\
Gaussian\vspace{0.3cm}
\end{minipage}\\
\begin{minipage}{0.32\linewidth}
\centering
\includegraphics[clip=true,height=1.7in]{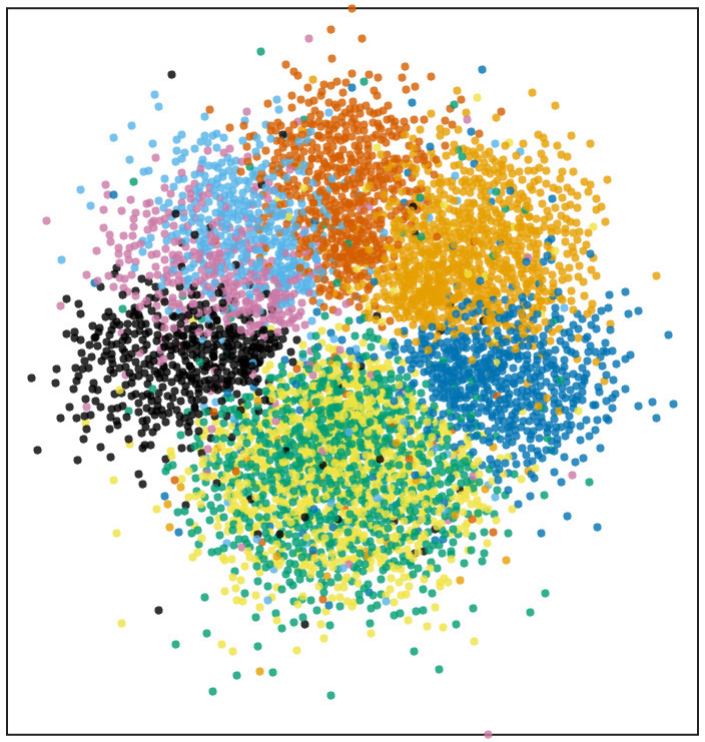}\\
CS
\end{minipage}
\begin{minipage}{0.32\linewidth}
\centering
\includegraphics[clip=true,height=1.7in]{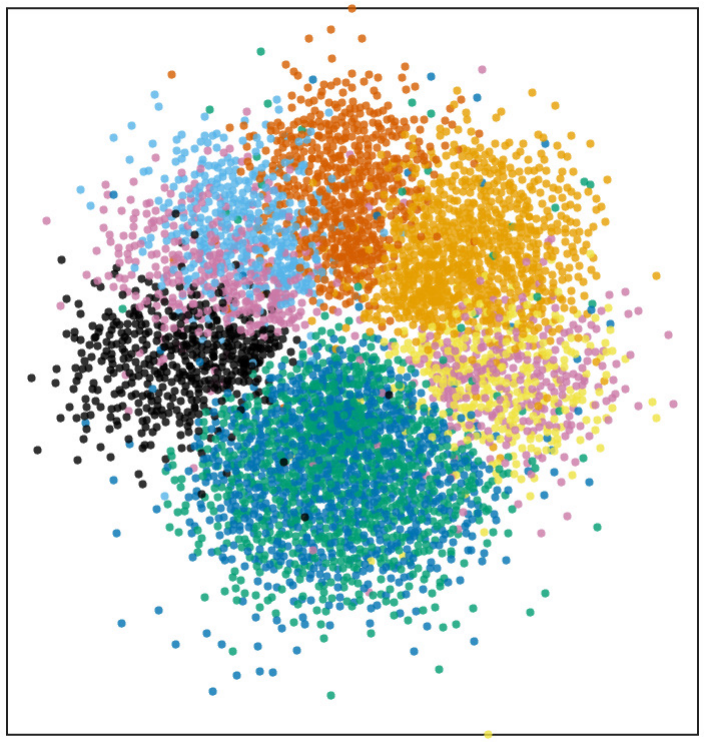}\\
DS-Gaussian
\end{minipage}
\begin{minipage}{0.32\linewidth}
\centering
\includegraphics[clip=true,height=1.7in]{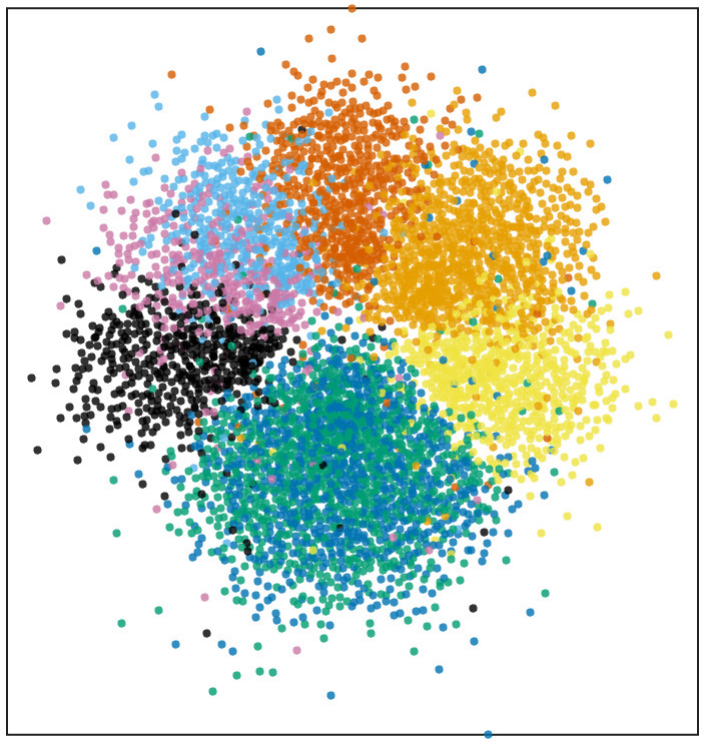}\\
DS-CS
\end{minipage}
\caption{Clustering results for the ABCD benchmark graph with $k=8$ clusters.}
\label{fig:synthetic-abcd}
\end{figure}
\begin{figure}[htbp]
\centering
\begin{minipage}{0.32\linewidth}
\centering
\includegraphics[clip=true,height=1.7in]{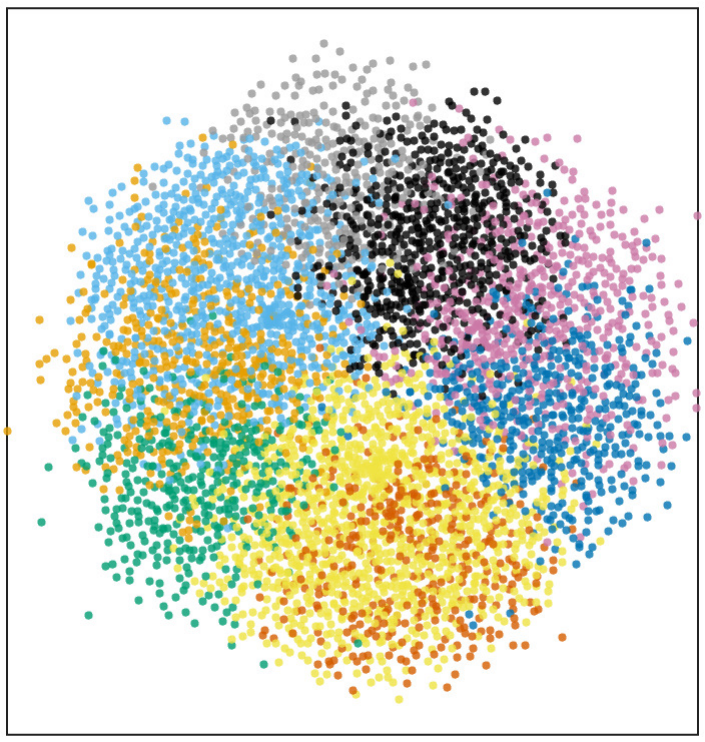}\\
Ground truth \vspace{0.3cm}
\end{minipage}
\begin{minipage}{0.32\linewidth}
\centering
\includegraphics[clip=true,height=1.7in]{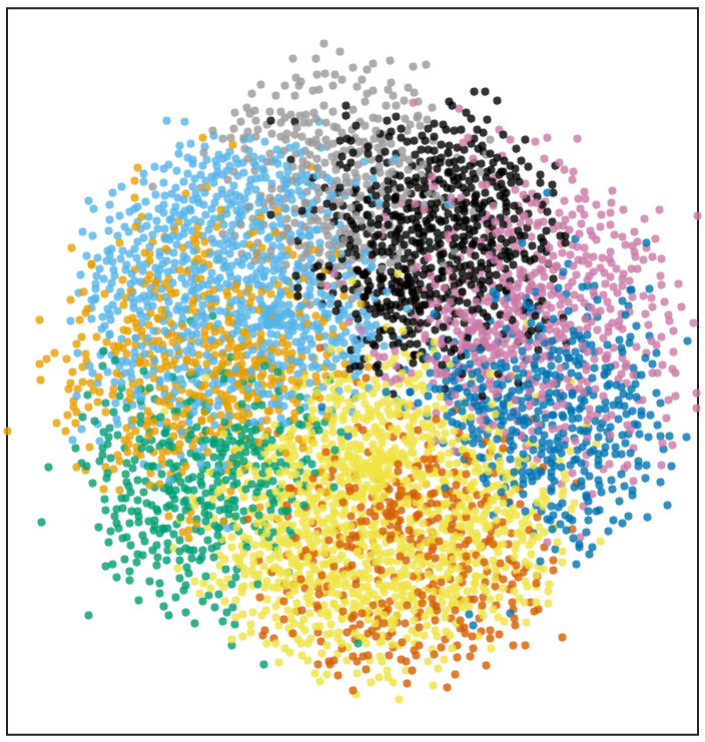}\\
Standard\vspace{0.3cm}
\end{minipage}
\begin{minipage}{0.32\linewidth}
\centering
\includegraphics[clip=true,height=1.7in]{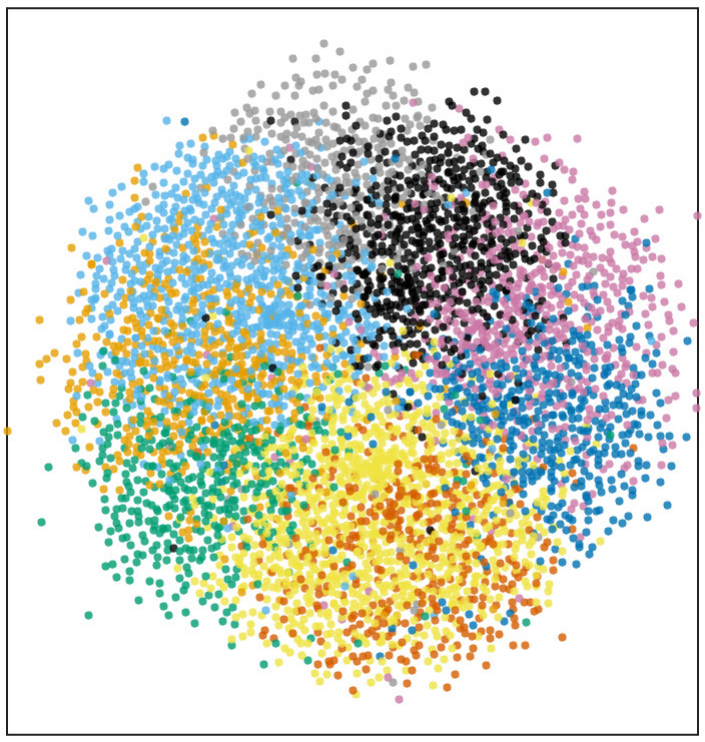}\\
Gaussian\vspace{0.3cm}
\end{minipage}\\
\begin{minipage}{0.32\linewidth}
\centering
\includegraphics[clip=true,height=1.7in]{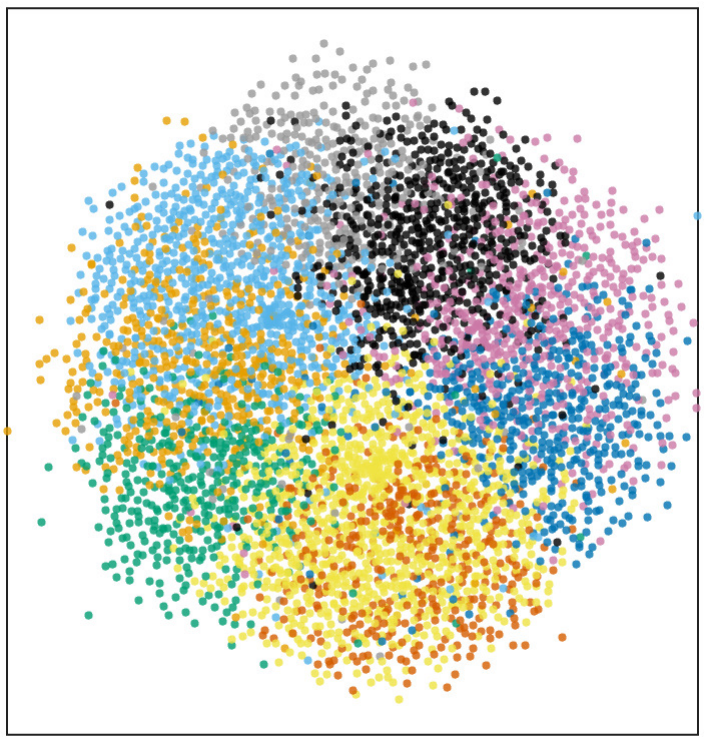}\\
CS
\end{minipage}
\begin{minipage}{0.32\linewidth}
\centering
\includegraphics[clip=true,height=1.7in]{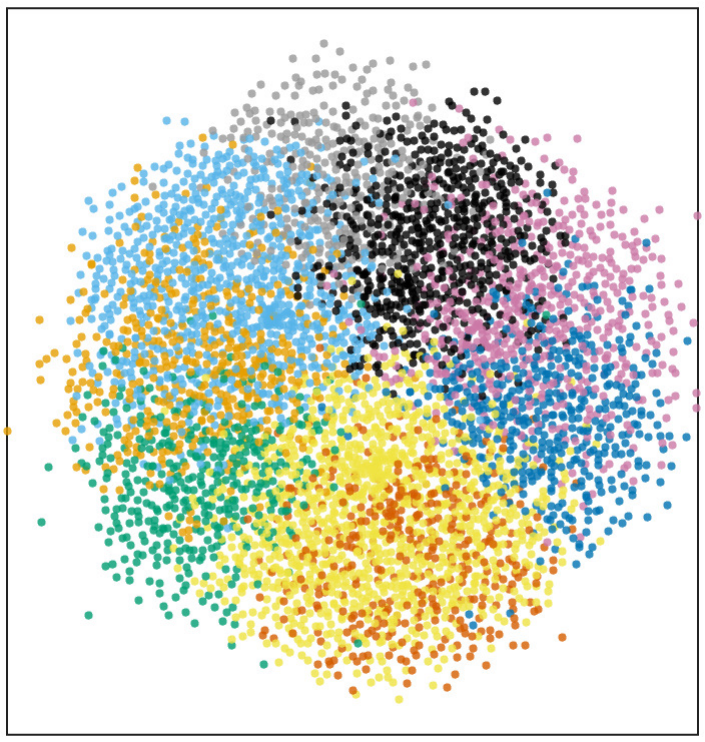}\\
DS-Gaussian
\end{minipage}
\begin{minipage}{0.32\linewidth}
\centering
\includegraphics[clip=true,height=1.7in]{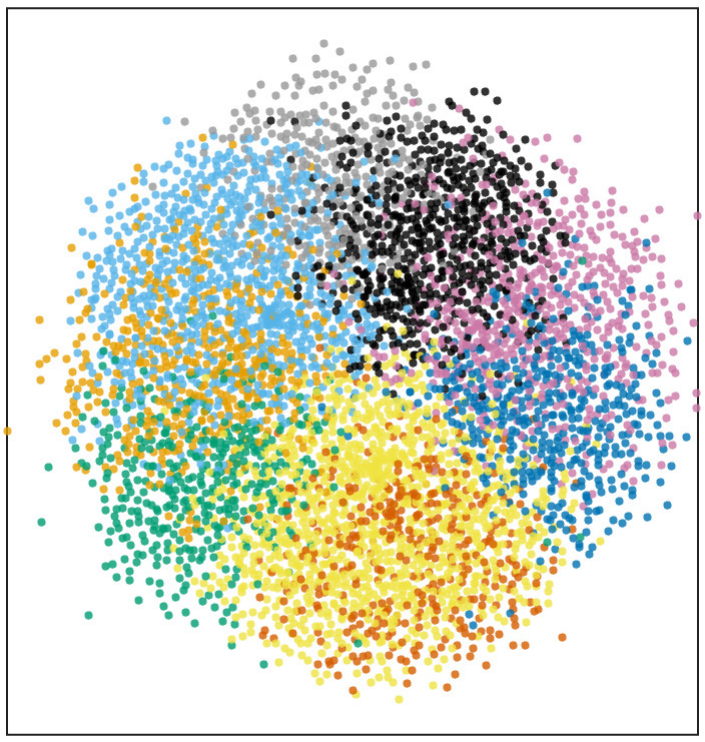}\\
DS-CS
\end{minipage}
\caption{Clustering results for the LFR benchmark graph with $k=9$ clusters.}
\label{fig:synthetic-lfr}
\end{figure}

To visualize the clustering results, we represent each graph as a scatter plot of its nodes, with node coordinates computed by the force-directed layout in MATLAB \cite{fruchterman1991graph} and fixed for all six panels. The panel labeled Ground truth is colored according to the ground-truth labels, and the other panels according to the predicted labels. For each method, we show the run whose Macro-F1 score is closest to the average over the $30$ runs, after permuting the predicted labels by the one-to-one matching that maximizes the sum of the matching scores weighted by the reciprocal of the ground-truth cluster sizes.
\Cref{fig:synthetic-dcsbm,fig:synthetic-abcd,fig:synthetic-lfr} show the clustering results for the synthetic DC-SBM, ABCD, and LFR benchmark graphs with $k=4$, $8$, and $9$ clusters, respectively. In each panel, the nodes with the same color belong to the same cluster after the label matching. For the synthetic DC-SBM and ABCD benchmark graphs, the DS-Gaussian and DS-CS panels resemble the Ground truth panel more closely than the Gaussian and CS panels do.

\subsection{Real-world benchmark}
To assess the randomized methods on real-world network data, we consider three graphs: Political Blogs, PubMed, and Facebook Page-Page.
The Political Blogs graph is constructed from the hyperlink network among U.S. political blogs \cite{adamic2005political}, where each node represents a blog and each edge a hyperlink between two blogs. We symmetrize the directed network to obtain an undirected graph and remove self-loops. Each node is labeled as liberal or conservative, and the labels are used as the ground-truth clusters. The largest connected component has $1222$ nodes, $\num{16714}$ edges, and $k=2$ ground-truth clusters, with a degree Gini coefficient of $0.6220$.
The PubMed graph is constructed from the PubMed citation network provided by the Planetoid dataset \cite{yang2016revisiting}. Its nodes are scientific publications on diabetes, and its edges are citation links between them. The network is symmetrized and self-loops are removed, as before. The three diabetes-related subject categories of the publications serve as the ground-truth clusters. In this case, the largest connected component contains $\num{19717}$ nodes and $\num{44324}$ edges, with $k=3$ ground-truth clusters and a degree Gini coefficient of $0.6037$.
The Facebook Page-Page graph is constructed from the network of verified Facebook pages in the MUSAE dataset \cite{rozemberczki2021multi}. Nodes correspond to Facebook pages and edges to mutual likes between them. We remove self-loops from the network. The four page categories define the ground-truth clusters. Restricting the graph to its largest connected component leaves $\num{22470}$ nodes and $\num{170823}$ edges, with $k=4$ ground-truth clusters and a degree Gini coefficient of $0.6164$.

\begin{figure}[htbp]
\centering
\begin{minipage}{0.32\linewidth}
\centering
\includegraphics[clip=true,height=1.7in]{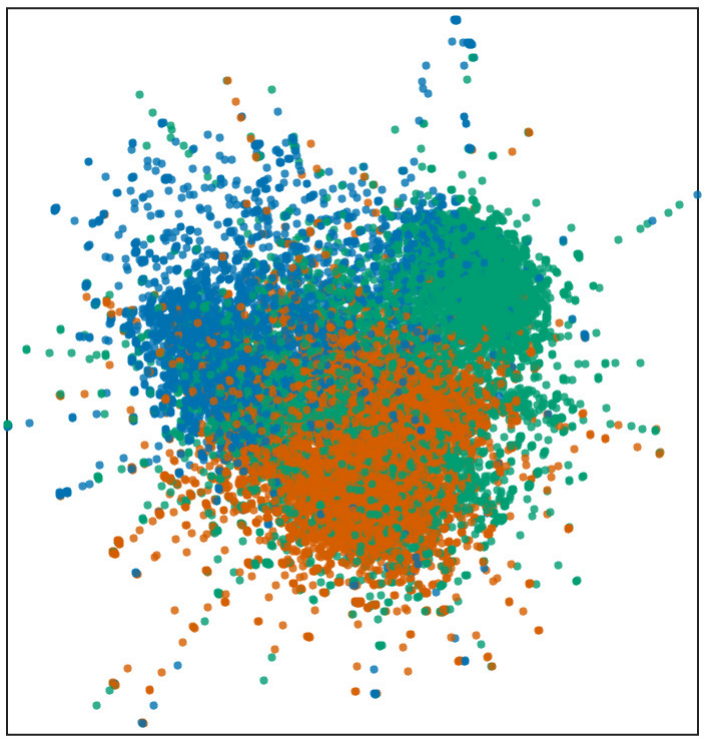}\\
Ground truth \vspace{0.3cm}
\end{minipage}
\begin{minipage}{0.32\linewidth}
\centering
\includegraphics[clip=true,height=1.7in]{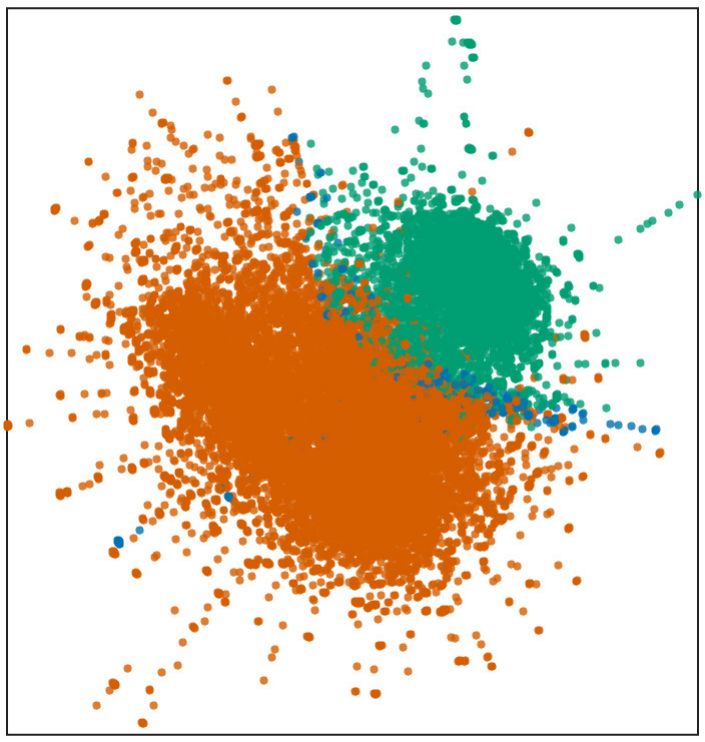}\\
Standard\vspace{0.3cm}
\end{minipage}
\begin{minipage}{0.32\linewidth}
\centering
\includegraphics[clip=true,height=1.7in]{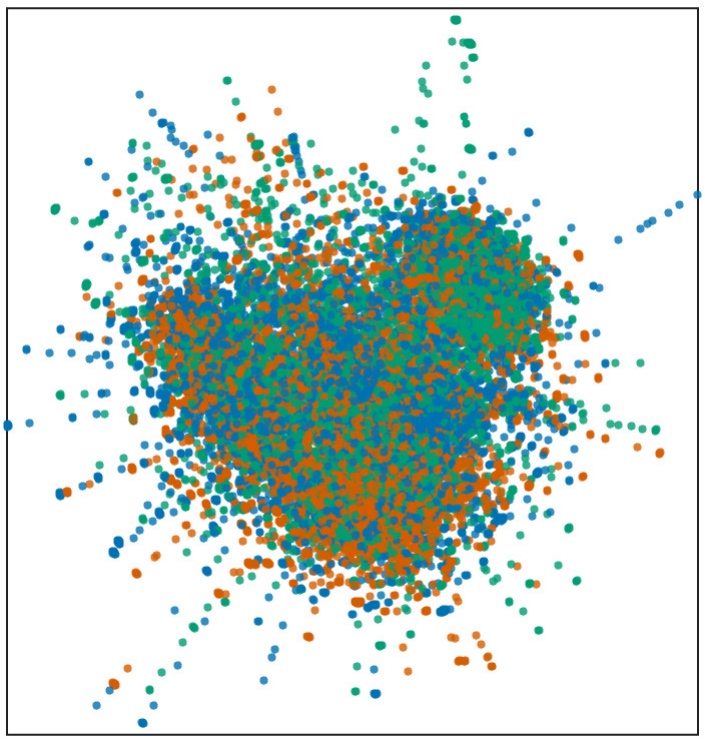}\\
Gaussian\vspace{0.3cm}
\end{minipage}\\
\begin{minipage}{0.32\linewidth}
\centering
\includegraphics[clip=true,height=1.7in]{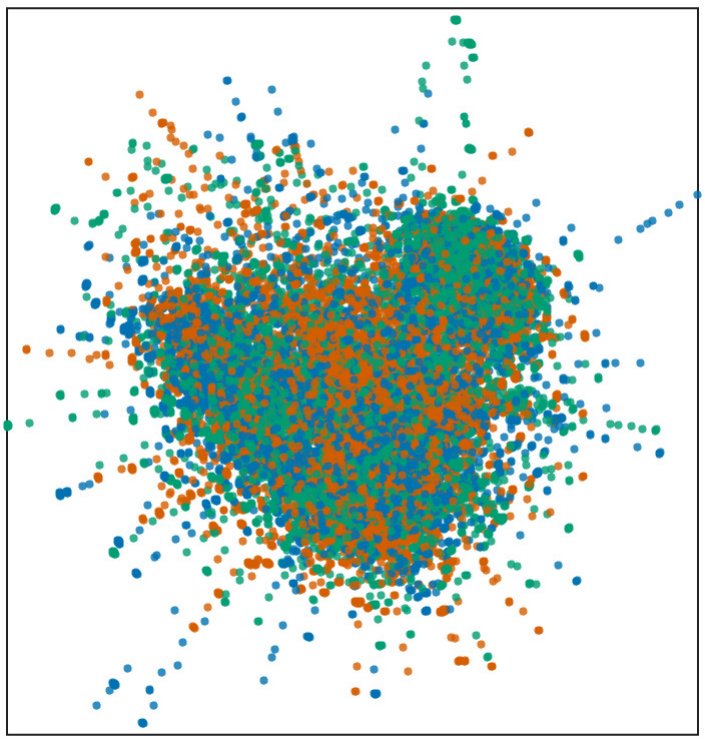}\\
CS
\end{minipage}
\begin{minipage}{0.32\linewidth}
\centering
\includegraphics[clip=true,height=1.7in]{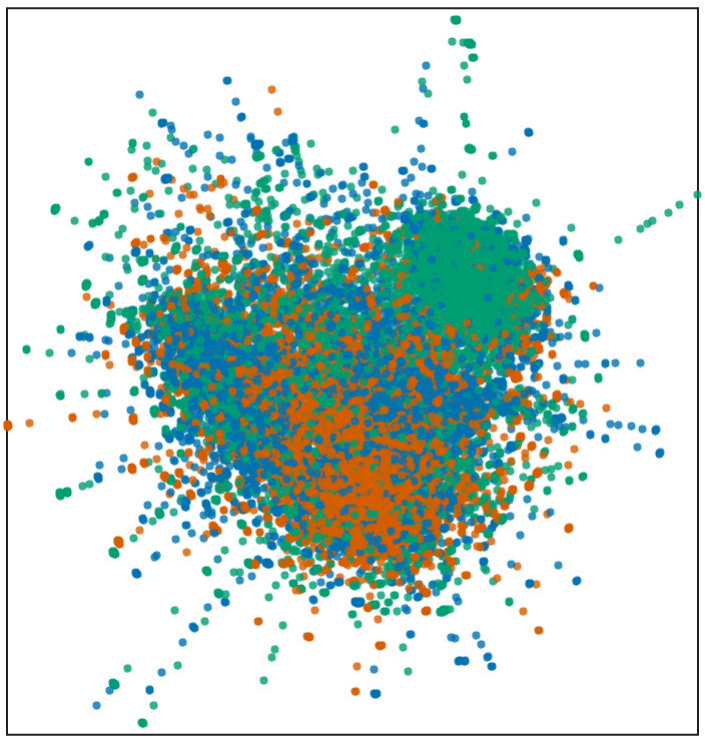}\\
DS-Gaussian
\end{minipage}
\begin{minipage}{0.32\linewidth}
\centering
\includegraphics[clip=true,height=1.7in]{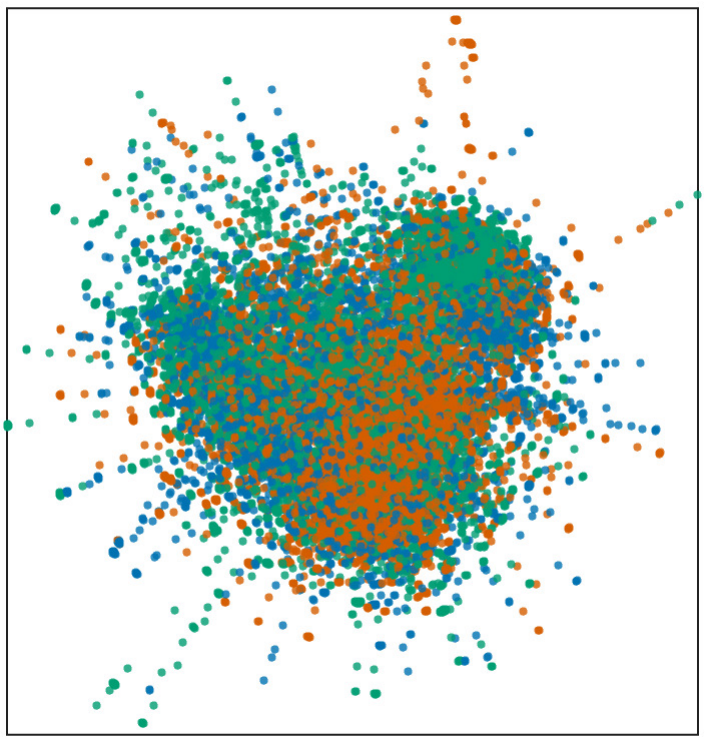}\\
DS-CS
\end{minipage}
\caption{Clustering results for the PubMed graph with $k=3$ clusters.}
\label{fig:realdata-pubmed}
\end{figure}

\begin{figure}[htbp]
\centering
\begin{minipage}{0.32\linewidth}
\centering
\includegraphics[clip=true,height=1.7in]{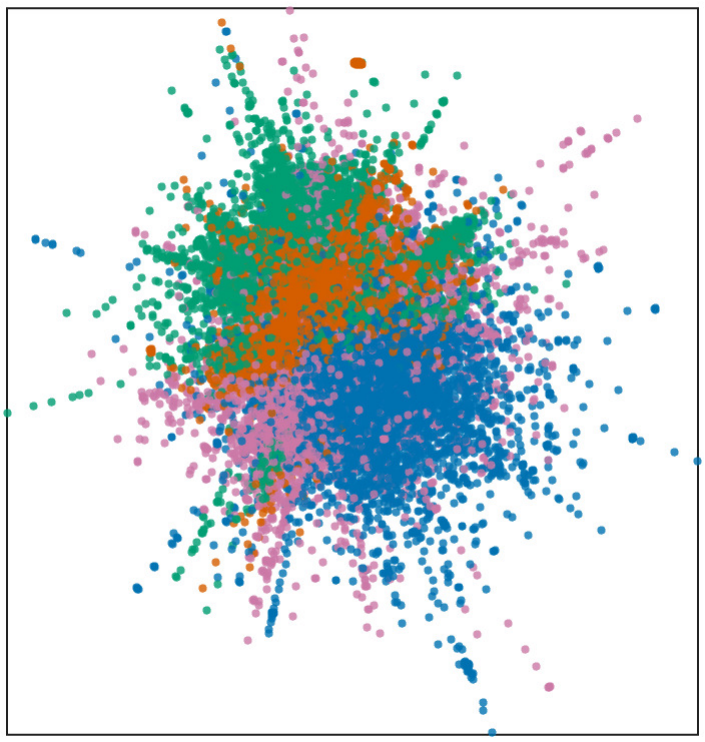}\\
Ground truth \vspace{0.3cm}
\end{minipage}
\begin{minipage}{0.32\linewidth}
\centering
\includegraphics[clip=true,height=1.7in]{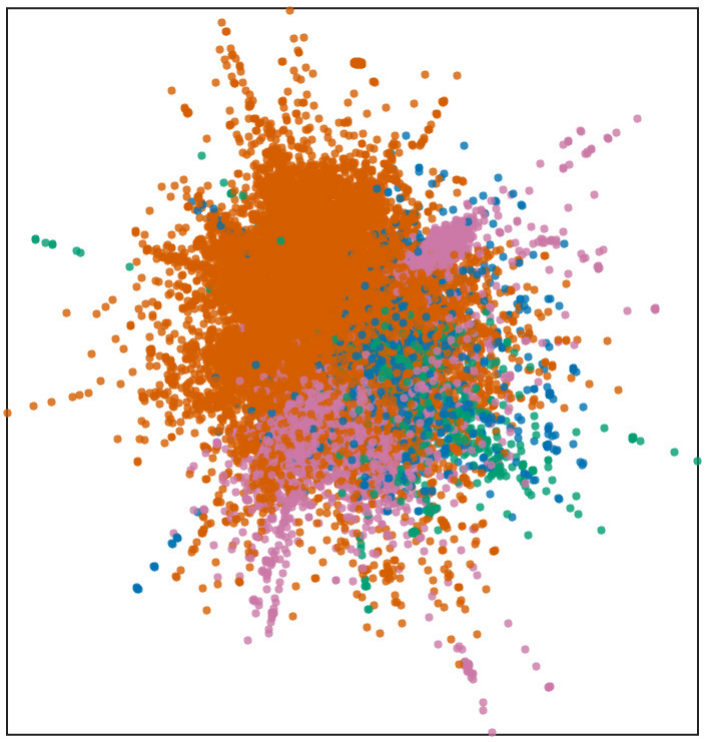}\\
Standard\vspace{0.3cm}
\end{minipage}
\begin{minipage}{0.32\linewidth}
\centering
\includegraphics[clip=true,height=1.7in]{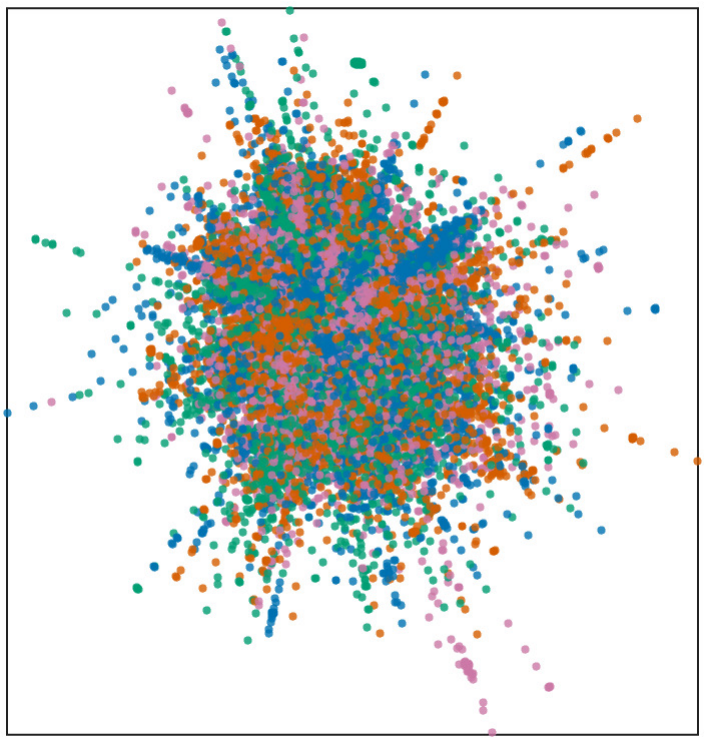}\\
Gaussian\vspace{0.3cm}
\end{minipage}\\
\begin{minipage}{0.32\linewidth}
\centering
\includegraphics[clip=true,height=1.7in]{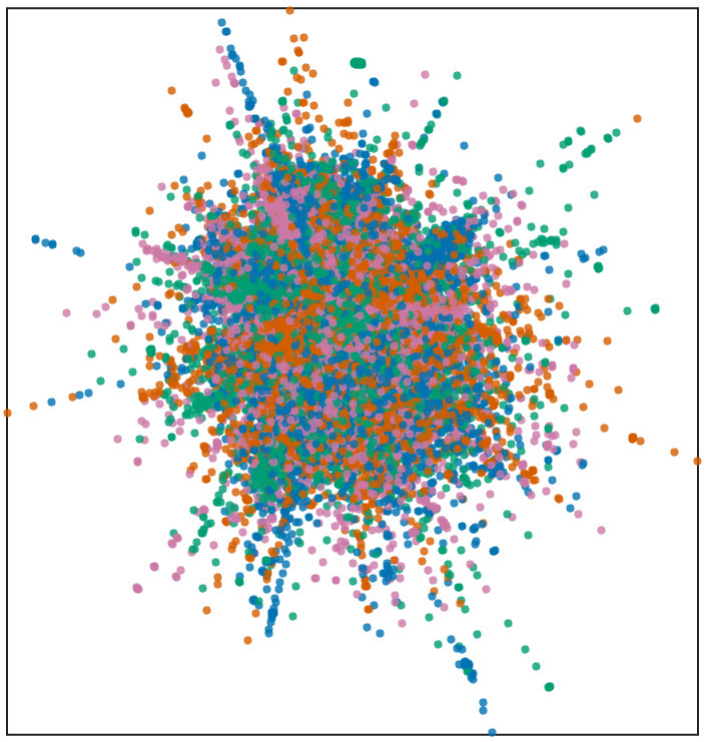}\\
CS
\end{minipage}
\begin{minipage}{0.32\linewidth}
\centering
\includegraphics[clip=true,height=1.7in]{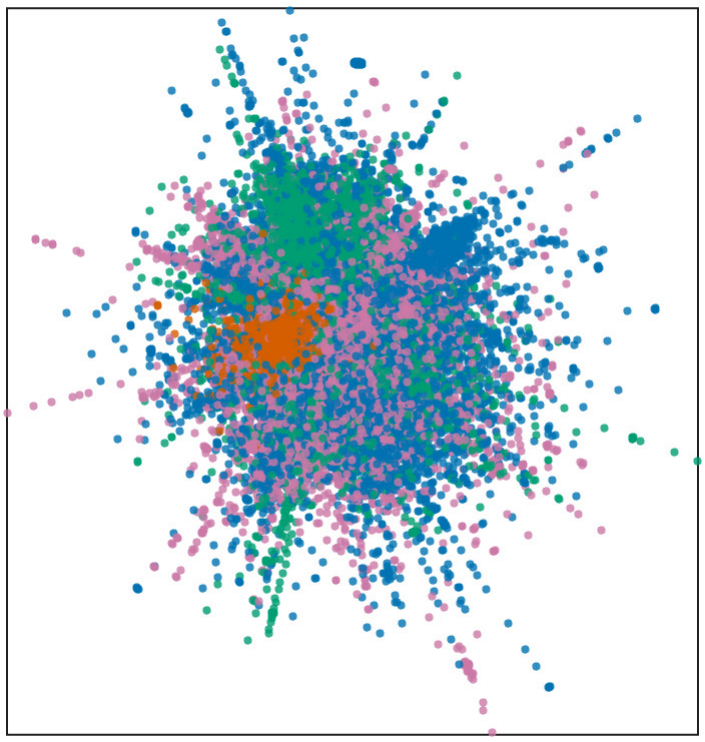}\\
DS-Gaussian
\end{minipage}
\begin{minipage}{0.32\linewidth}
\centering
\includegraphics[clip=true,height=1.7in]{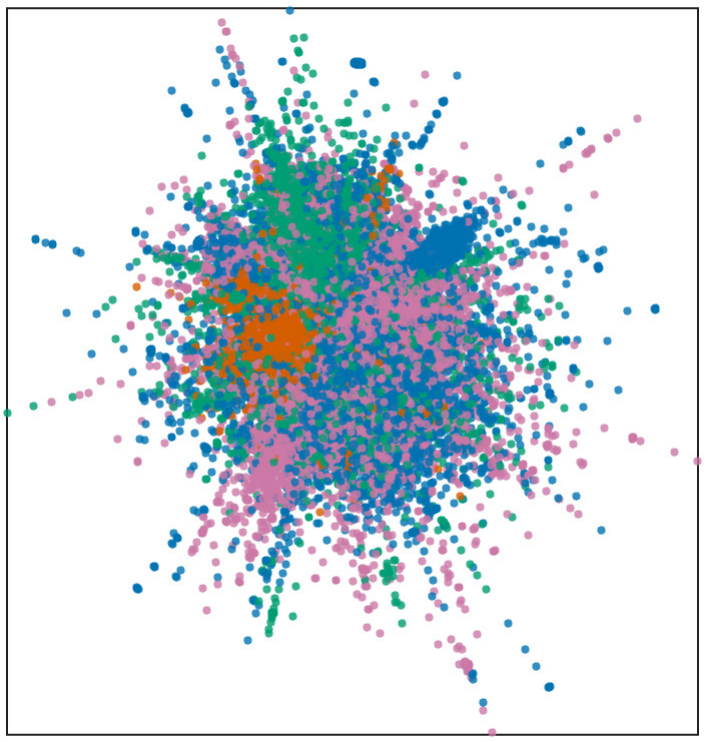}\\
DS-CS
\end{minipage}
\caption{Clustering results for the Facebook Page-Page graph with $k=4$ clusters.}
\label{fig:realdata-facebook}
\end{figure}

\Cref{tab:accuracy-realdata} reports the corresponding results for the three real-world graphs. The reference values of Macro-F1, NMI, and ARI obtained by standard spectral clustering are $0.9573$, $0.7492$, and $0.8369$ for Political Blogs, $0.4437$, $0.1740$, and $0.1432$ for PubMed, and $0.3017$, $0.1658$, and $0.0988$ for Facebook Page-Page, respectively.
On all three graphs, the DS methods are more accurate than the conventional Gaussian and CS methods in every measure. For Political Blogs, DS-Gaussian attains between $96.70\%$ and $99.40\%$ of the reference values. For PubMed and Facebook Page-Page, the NMI and ARI values of all five methods are low, which indicates that the ground-truth labels of these graphs are only weakly aligned with the leading eigenvectors of $T_\perp$. The NMI and ARI values of DS-Gaussian amount to $29.89\%$ and $40.50\%$ of the reference for PubMed and to $59.17\%$ and $72.57\%$ for Facebook Page-Page, whereas those of the conventional Gaussian and CS methods remain close to zero. The Macro-F1 scores of DS-Gaussian on these two graphs are $101.40\%$ and $132.81\%$ of the reference values. Since standard spectral clustering assigns most nodes to a single cluster here (see \Cref{fig:realdata-pubmed,fig:realdata-facebook}), this should not be read as a better recovery of the exact spectral embedding. Overall, with the same sketch dimension, the DS methods recover substantially more of the accuracy of standard spectral clustering than their conventional counterparts.

\begin{table}[htbp]
\centering
\footnotesize
\caption{Average clustering accuracy over $30$ runs on the real-world graphs.}
\label{tab:accuracy-realdata}
\begin{tabular}{cc cccc}
\toprule
Graph & Metric & Gaussian & CS & DS-Gaussian & DS-CS \\
\midrule
\multirow{3}{*}{Political Blogs}
& Macro-F1  & 0.9429 & 0.9419 & {\bf 0.9516} (99.40\%) & 0.9504 \\
& NMI & 0.6896 & 0.6859 & {\bf 0.7245} (96.70\%) & 0.7192 \\
& ARI & 0.7852 & 0.7819 & {\bf 0.8162} (97.53\%) & 0.8120 \\
\midrule
\multirow{3}{*}{PubMed}
& Macro-F1  & 0.3653 & 0.3780 & {\bf 0.4499} (101.40\%) & 0.4429 \\
& NMI & 0.0068 & 0.0101 & {\bf 0.0520} (29.89\%) & 0.0461 \\
& ARI & 0.0079 & 0.0116 & {\bf 0.0580} (40.50\%) & 0.0539 \\
\midrule
\multirow{3}{*}{Facebook Page-Page}
& Macro-F1  & 0.2949 & 0.2991 & {\bf 0.4007} (132.81\%) & 0.3973 \\
& NMI & 0.0117 & 0.0140 & {\bf 0.0981} (59.17\%) & 0.0935 \\
& ARI & 0.0114 & 0.0135 & {\bf 0.0717} (72.57\%) & 0.0700 \\
\bottomrule
\end{tabular}
\end{table}

\begin{figure}[htbp]
\centering
\begin{minipage}{0.32\linewidth}
\centering
\includegraphics[clip=true,height=1.7in]{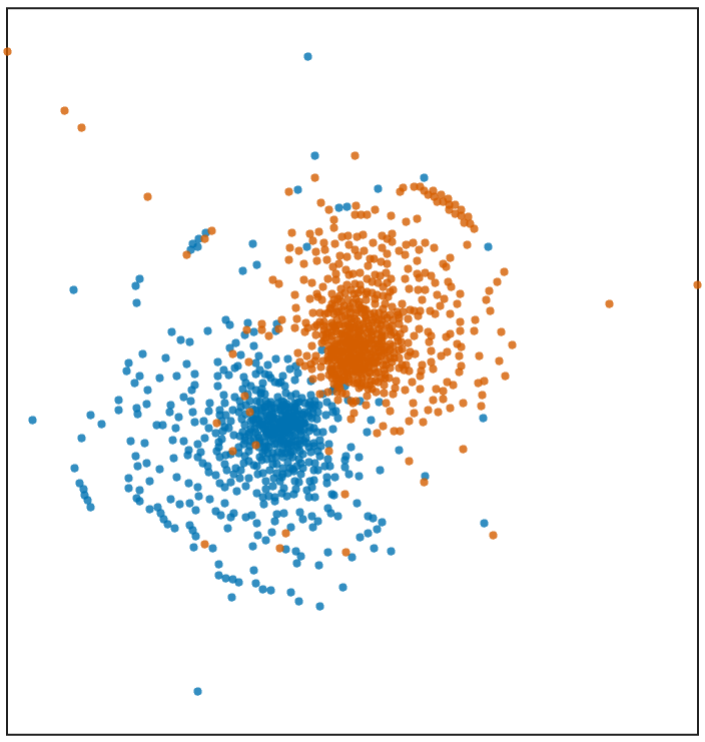}\\
Ground truth \vspace{0.3cm}
\end{minipage}
\begin{minipage}{0.32\linewidth}
\centering
\includegraphics[clip=true,height=1.7in]{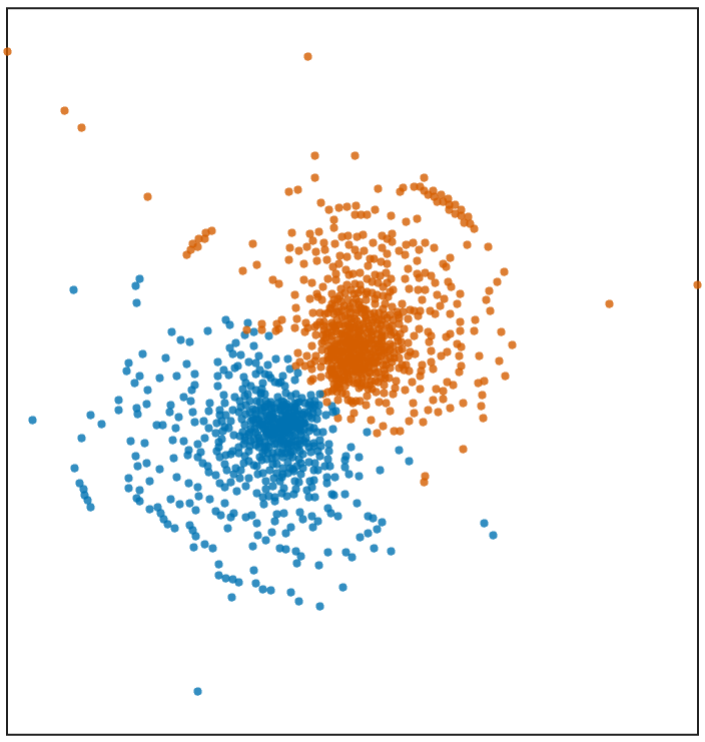}\\
Standard\vspace{0.3cm}
\end{minipage}
\begin{minipage}{0.32\linewidth}
\centering
\includegraphics[clip=true,height=1.7in]{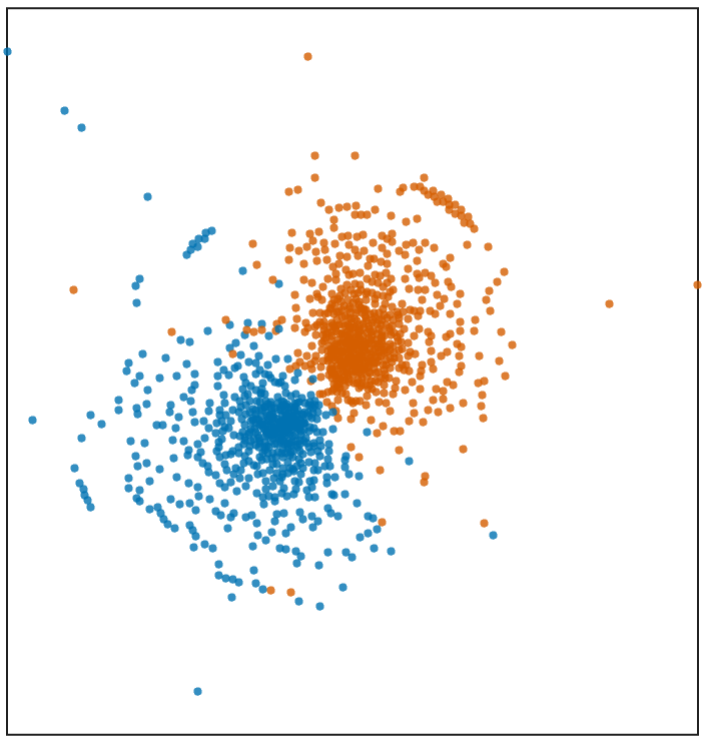}\\
Gaussian\vspace{0.3cm}
\end{minipage}\\
\begin{minipage}{0.32\linewidth}
\centering
\includegraphics[clip=true,height=1.7in]{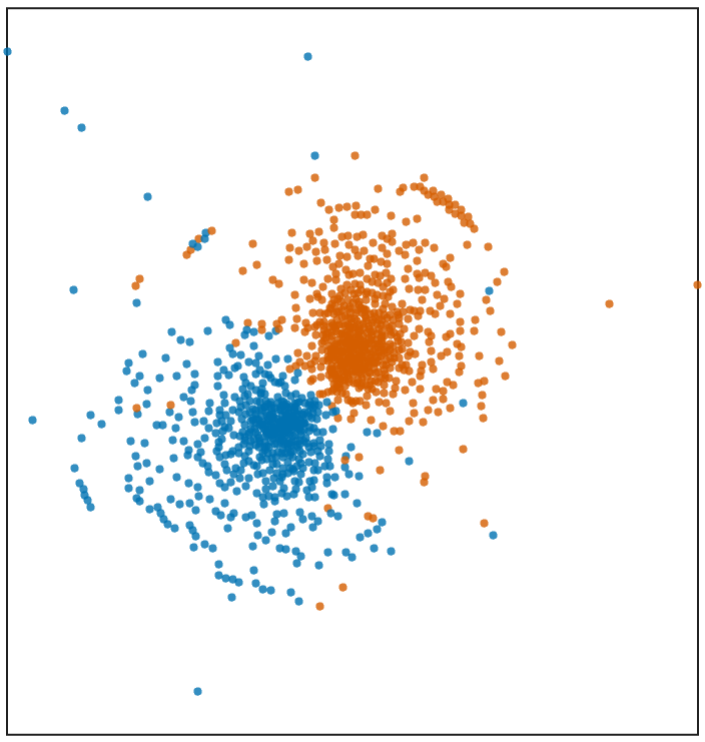}\\
CS
\end{minipage}
\begin{minipage}{0.32\linewidth}
\centering
\includegraphics[clip=true,height=1.7in]{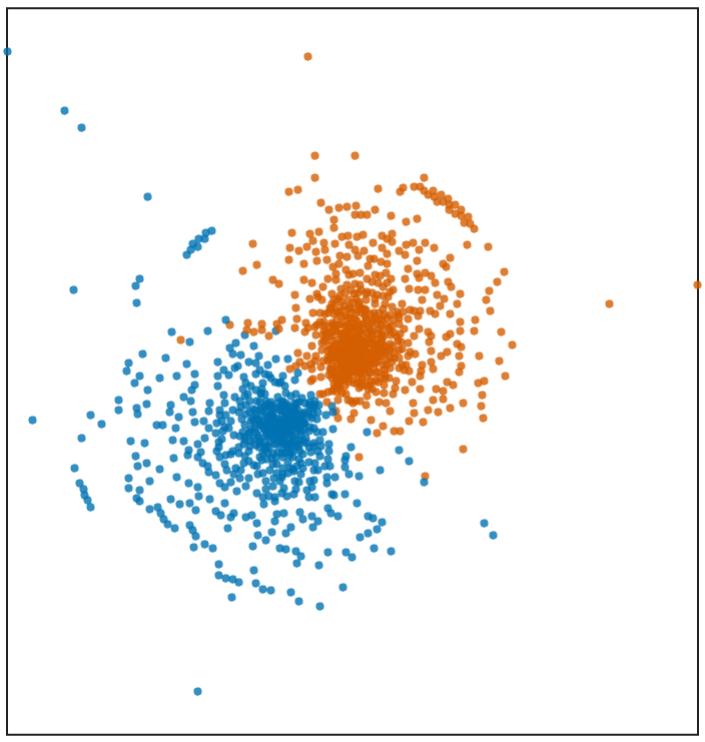}\\
DS-Gaussian
\end{minipage}
\begin{minipage}{0.32\linewidth}
\centering
\includegraphics[clip=true,height=1.7in]{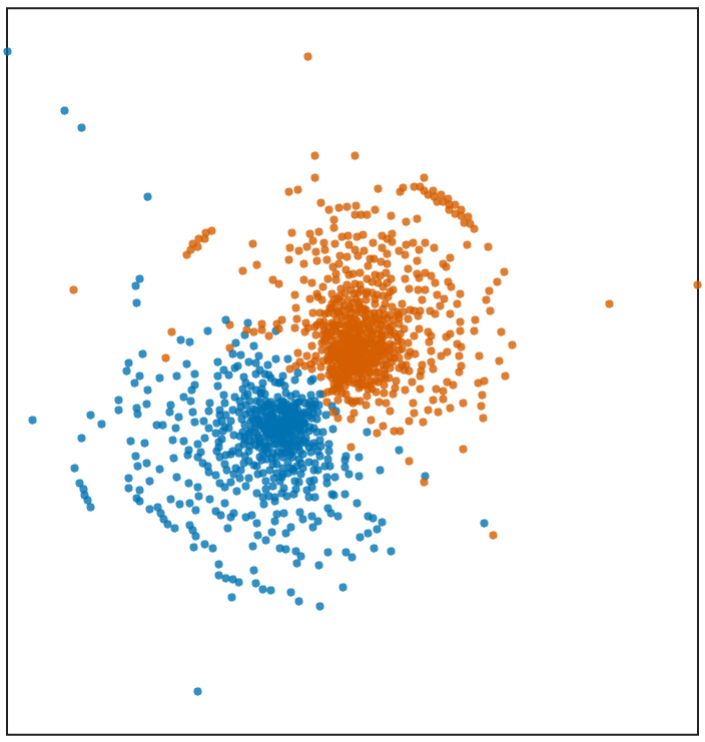}\\
DS-CS
\end{minipage}
\caption{Clustering results for the Political Blogs graph with $k=2$ clusters.}
\label{fig:realdata-politicalblogs}
\end{figure}


We visualize the clustering results for the real-world graphs in the same manner as for the synthetic benchmark graphs. 
\Cref{fig:realdata-politicalblogs,fig:realdata-pubmed,fig:realdata-facebook} display the Political Blogs, PubMed, and Facebook Page-Page graphs with $k=2$, $3$, and $4$ clusters, respectively. For the Political Blogs graph, all five methods give results similar to the ground truth. For the PubMed and Facebook Page-Page graphs, standard spectral clustering assigns most of the nodes to one cluster and hardly identifies one of the ground-truth clusters. For the PubMed graph, the DS-Gaussian and DS-CS panels are closer to the Ground truth panel than the Standard, Gaussian, and CS panels.

\subsection{Large-scale experiments}

We finally compare the computational efficiency of the five methods on two large-scale real-world graphs, Pokec and LiveJournal, from the SNAP collection \cite{leskovec2014snap}. The graphs are constructed from the directed networks of \cite{takac2012data} and \cite{backstrom2006group,leskovec2009community}, respectively. We convert each directed network into an undirected graph, replace multiple edges between the same pair of nodes with a single edge of unit weight, and remove self-loops. The resulting Pokec graph is connected and has $\num{1632803}$ nodes and $\num{22301964}$ edges. For LiveJournal, we consider the largest connected component, which has $\num{4843953}$ nodes and $\num{42845684}$ edges.

For both graphs, we set $k=4$. We use $q=2$ and $r=20$ for Pokec and $q=4$ and $r=30$ for LiveJournal. Thus, the sketch dimensions are $\ell=24$ and $\ell=34$, respectively. For standard spectral clustering, we compute the $k$ leading eigenvectors of $T_\perp$ using ARPACK \cite{lehoucq1998arpack}, where the maximum dimension of the Krylov subspace and the convergence tolerance are set to $40$ and $10^{-6}$, respectively. All computations are performed using a single thread, and each experiment is repeated three times. The embedding times and Ritz values are averaged over the three runs.

The embedding time of standard spectral clustering includes the computation of the $k$ leading eigenvectors of $T_\perp$ and the application of $P_\perp$ to the computed eigenvectors. For the randomized methods, it comprises the construction of the degree strata for the DS methods, the generation of the sketching matrices, the computation of $Y$, the subspace iteration, and the Rayleigh--Ritz projection. The values in parentheses in the tables are the ratios of the embedding time of standard spectral clustering to those of the randomized methods.

Since ground-truth labels are not used for these graphs, we assess the accuracy through the Ritz values $\theta_1\geq\theta_2\geq\cdots$ of $C=Q^\top T_\perp Q$. By the Cauchy interlacing theorem, $\theta_i\leq\mu_{i+1}$ for $i=1,\ldots,\ell$ \cite{parlett1998symmetric}, so that a larger Ritz value is necessarily closer to the corresponding eigenvalue of $T_\perp$. For comparison, the four randomized methods are also applied to the adjacency matrix $A$ with the same values of $q$ and $r$, where $\lambda_i(A)$ denotes the $i$th largest eigenvalue of $A$. The $20$ largest eigenvalues of $A$ and $T_\perp$ are computed using ARPACK, with the maximum Krylov subspace dimension and the convergence tolerance set to $100$ and $10^{-6}$, respectively.

\begin{figure}[htbp]
\centering
\begin{minipage}{0.5\linewidth}
\centering
\includegraphics[clip=true,height=1.8in]{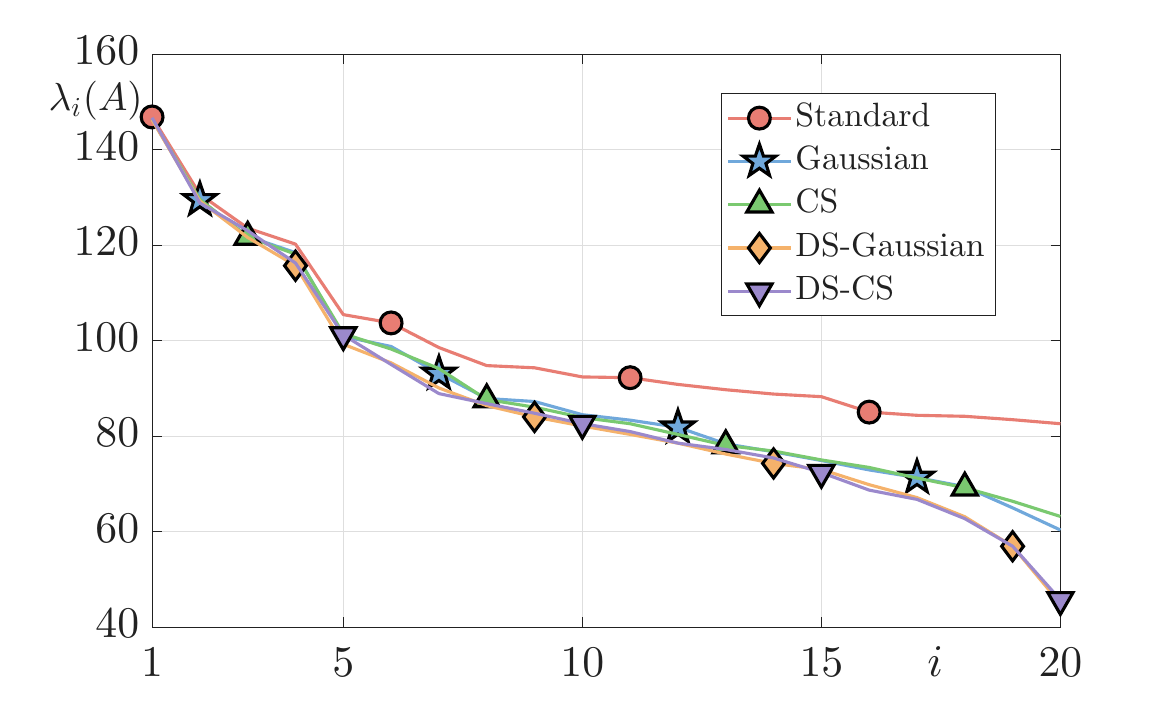}
\end{minipage}
\begin{minipage}{0.49\linewidth}
\centering
\includegraphics[clip=true,height=1.8in]{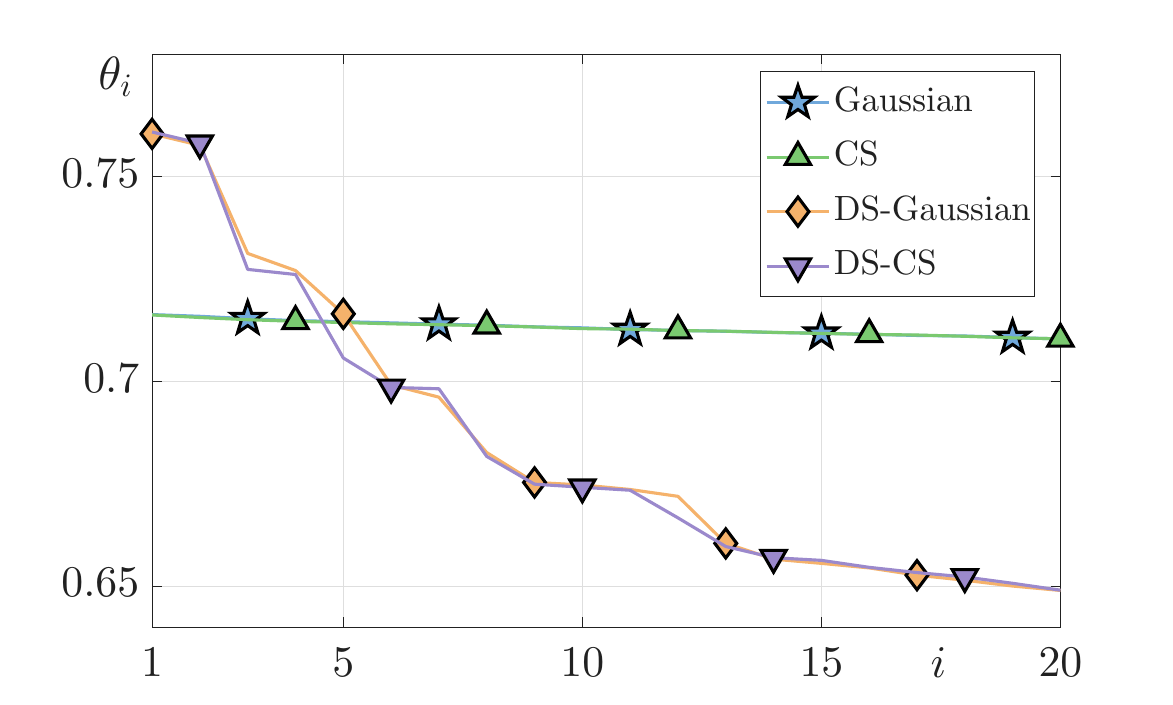}
\end{minipage}
\caption{Largest $20$ eigenvalues and Ritz values of $A$ (left) and largest $20$ Ritz values of $T_\perp$ (right) of the four randomized methods for the Pokec graph.}
\label{fig:spectrum-pokec}
\end{figure}

For the Pokec graph, \Cref{tab:efficiency-pokec} shows that all four randomized methods are faster than standard spectral clustering, and \Cref{fig:spectrum-pokec} compares the leading eigenvalues and Ritz values. For $A$, they all reproduce the overall trend of the exact spectrum. For $T_\perp$, Gaussian and CS yield nearly flat Ritz spectra, which indicates that, with the present values of $q$ and $\ell$, the sampled range is dominated by the bulk of the spectrum. By contrast, the leading Ritz values of DS-Gaussian and DS-CS, in particular those with $i\leq k$, are larger than those of Gaussian and CS. In view of the interlacing property, the DS methods approximate the dominant eigenvalues of $T_\perp$ more accurately with the same sketch dimension, in agreement with the analysis in \Cref{sec:analysis}.

\begin{table}[htbp]
\centering
\footnotesize
\caption{Average embedding times (speed up) of the five methods for the Pokec graph.}
\label{tab:efficiency-pokec}
\begin{tabular}{cccccc}
\toprule
 & Standard & Gaussian & CS & DS-Gaussian & DS-CS \\
\midrule
Embedding time (s) & 21.928 & 13.045 (1.68) & 12.565 (1.75) & 13.827 (1.59) & 12.540 (1.75) \\
\bottomrule
\end{tabular}
\end{table}

\begin{figure}[htbp]
\centering
\begin{minipage}{0.5\linewidth}
\centering
\includegraphics[clip=true,height=1.8in]{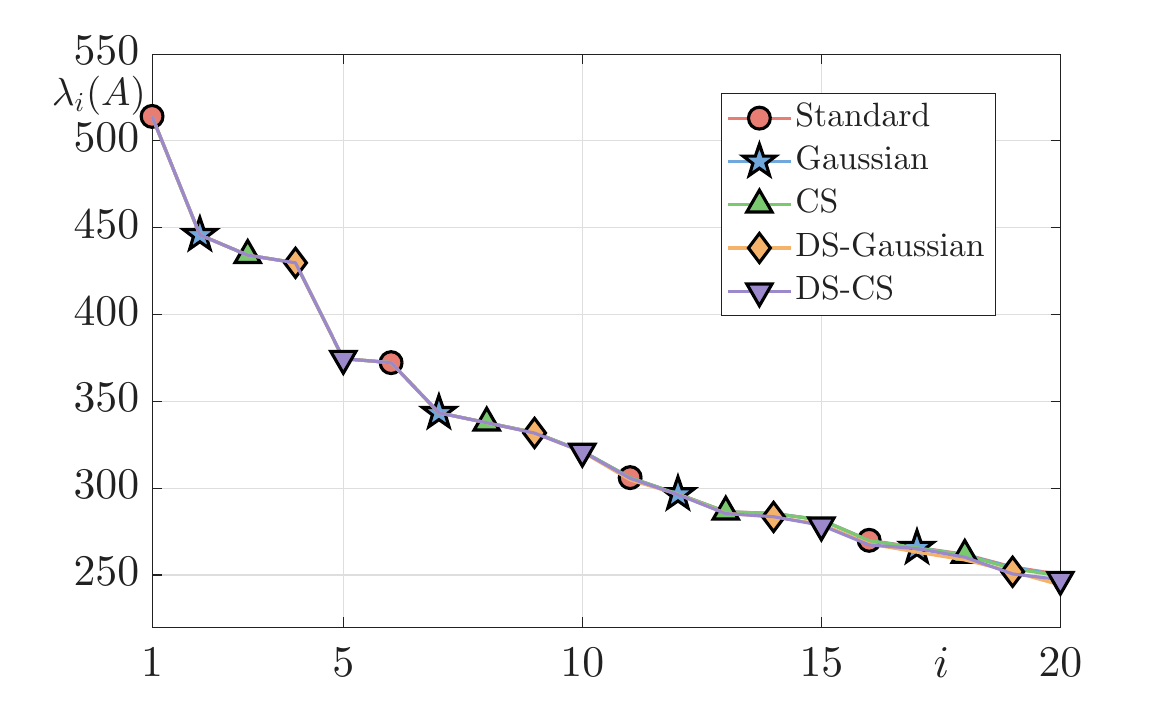}
\end{minipage}
\begin{minipage}{0.49\linewidth}
\centering
\includegraphics[clip=true,height=1.8in]{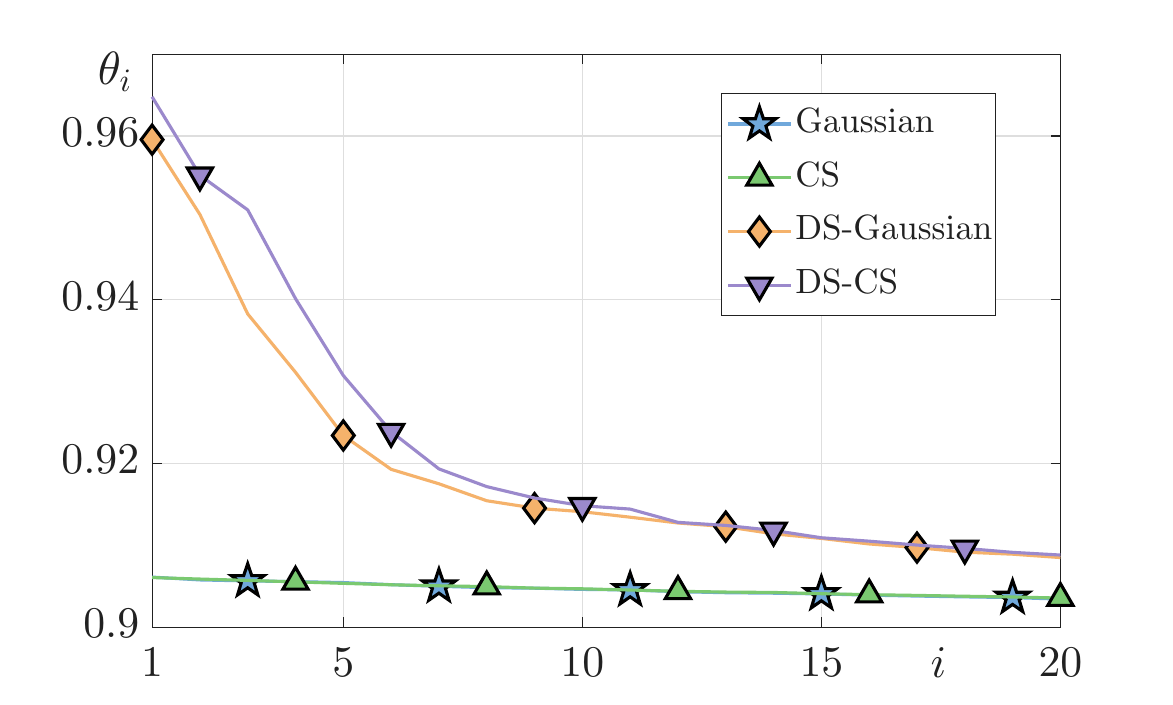}
\end{minipage}
\caption{Largest $20$ eigenvalues and Ritz values of $A$ (left) and largest $20$ Ritz values of $T_\perp$ (right) of the four randomized methods for the LiveJournal graph.}
\label{fig:spectrum-livejournal}
\end{figure}

For LiveJournal, we observe the same qualitative behavior. As reported in \Cref{tab:efficiency-livejournal}, all four randomized methods reduce the embedding time by approximately a factor of eight compared with standard spectral clustering. In \Cref{fig:spectrum-livejournal}, Gaussian and CS again yield nearly flat Ritz spectra for $T_\perp$, whereas the Ritz values of DS-Gaussian and DS-CS are larger over the whole range of indices shown. On both large-scale graphs, degree stratification thus yields a sampled range that captures the leading eigenspace of $T_\perp$ more accurately for a fixed sketch dimension and a fixed number of subspace iterations.

\begin{table}[htbp]
\centering
\footnotesize
\caption{Average embedding times (speed up) of the five methods for the LiveJournal graph.}
\label{tab:efficiency-livejournal}
\begin{tabular}{cccccc}
\toprule
 & Standard & Gaussian & CS & DS-Gaussian & DS-CS \\
\midrule
Embedding time (s) & 605.688 & 79.150 (7.65) & 76.662 (7.90) & 80.084 (7.56) & 76.490 (7.92) \\
\bottomrule
\end{tabular}
\end{table}

\section{Conclusion}\label{sec:conclusion}

Degree heterogeneity is an important structural feature of many real-world networks, yet it is not used in conventional randomized spectral clustering. We proposed degree-stratified randomized projection, which partitions the nodes into degree strata and distributes a fixed sketch dimension across the strata. Our analysis shows that the allocation acts through the leverage that each stratum carries for the target eigenspace, that this leverage is proportional to the degree volume under the degree-corrected stochastic block model, and that the accuracy of the sampled range carries over to the Ritz vectors and to the spectral embedding. In the numerical experiments, degree stratification improved the clustering accuracy of conventional randomized projection methods at the same sketch dimension and essentially the same computational cost.

Several questions remain open. The square-root allocation rule is a heuristic choice, and the optimal allocation is an experimental design problem that depends on the unknown stratum leverage matrices. Adaptive rules that estimate these quantities during the iteration are a natural next step. An analysis of the CountSketch variant and of structure-aware sketches that exploit graph characteristics beyond the node degrees is left for future work.

%%%%%%%%%%%%%%%%%%%%%%%%%%%%%%%%%%%%%%%%%%%%%%%%%%
% Appendices (the supplementary materials of the journal version)
\appendix

%%%%%%%%%%%%%%%%%%%%%%%%%%%%%%%%%%%%%%%%%%%%%%%%%%
\section{Proof for Section~\ref{sec:method}}
\label{supp:method}

\begin{proof}[Proof of \cref{prop:shift-deflation}]
Since $T=\frac12(I+S)$, for each eigenpair $S\bs{u}_i=\lambda_i\bs{u}_i$ we have
\[
    T\bs{u}_i=\frac12(I+S)\bs{u}_i=\frac12(1+\lambda_i)\bs{u}_i=\mu_i\bs{u}_i .
\]
Since $\lambda_i\in[-1,1]$ and $\lambda_1=1$, it follows that $\mu_i\in[0,1]$ and $\mu_1=1$. The map $\lambda\mapsto\frac12(1+\lambda)$ is strictly increasing, so it preserves the ordering of the eigenvalues.

For the deflated operator, $P_\perp\bs{u}_1=0$ gives $T_\perp\bs{u}_1=P_\perp TP_\perp\bs{u}_1=0$. For $i\geq 2$, orthogonality gives $P_\perp\bs{u}_i=\bs{u}_i$, and hence
\[
    T_\perp\bs{u}_i=P_\perp T\bs{u}_i=\mu_iP_\perp\bs{u}_i=\mu_i\bs{u}_i .
\]
\end{proof}

%%%%%%%%%%%%%%%%%%%%%%%%%%%%%%%%%%%%%%%%%%%%%%%%%%
\section{Proofs for Section~\ref{sec:analysis}}
\label{supp:analysis}

Throughout this section, let $Z_c=[\bs{u}_{k_\star+2},\ldots,\bs{u}_n,\bs{u}_1]$ contain the eigenvectors of $T_\perp$ that are not in $Z_\star$, and let $\Sigma_\star=\diag(\mu_2,\ldots,\mu_{k_\star+1})$ and $\Sigma_c=\diag(\mu_{k_\star+2},\ldots,\mu_n,0)$, so that
\[
    T_\perp=Z_\star\Sigma_\star Z_\star^\top+Z_c\Sigma_cZ_c^\top,
\]
where the matrix $[Z_\star,Z_c]$ is orthogonal. We write $\Pi=QQ^\top$ for the orthogonal projector onto $\range(Q)$, so that $\dist(Z_\star,Q)=\norm{(I-\Pi)Z_\star}_2$, and $\sigma_{\min}(\cdot)$ for the smallest singular value of a matrix. Since orthonormalization does not change the column space, the matrix $Q$ computed by \cref{alg:degree-stratified-rp} satisfies $\range(Q)=\range(T_\perp^{2q+1}\Omega)$.

To avoid reordering the nodes, we let $P_j\in\mathbb R^{n\times n_j}$ denote the column-selection matrix of stratum $\mathcal B_j$, so that $(T_\perp)_{:,\mathcal B_j}=T_\perp P_j$ and the test matrix \eqref{eq:structured-test-matrix} is
\[
    \Omega=[P_1\Omega_1,\ldots,P_s\Omega_s].
\]
The matrix $V_j=P_j^\top Z_\star\in\mathbb R^{n_j\times k_\star}$ contains the rows of $Z_\star$ that belong to $\mathcal B_j$; hence $W_j=V_j^\top V_j$, and we set $k_j=\rank(V_j)\leq k_\star$.

\subsection{Structural bound}

The following bound is a variant of the structural bounds for randomized subspace iteration \cite{halko2011finding,saibaba2019randomized}. It holds for an arbitrary test matrix and makes explicit that the columns of $\Omega$ may be rescaled freely.

\begin{lemma}
\label{lem:structural}
Let $q\geq 0$, let $\Omega\in\mathbb R^{n\times\ell}$ be a matrix for which $Z_\star^\top\Omega$ has full row rank, and let $Q$ be an orthonormal basis of $\range(T_\perp^{2q+1}\Omega)$. Then, for every nonsingular matrix $\Delta\in\mathbb R^{\ell\times\ell}$,
\[
    \dist(Z_\star,Q)
    \leq
    \left(\frac{\mu_{k_\star+2}}{\mu_{k_\star+1}}\right)^{2q+1}
    \frac{\norm{Z_c^\top\Omega\Delta}_2}{\sigma_{\min}(Z_\star^\top\Omega\Delta)} .
\]
\end{lemma}

\begin{proof}
Let $p=2q+1$, $\Omega_\star=Z_\star^\top\Omega\Delta\in\mathbb R^{k_\star\times\ell}$, and $\Omega_c=Z_c^\top\Omega\Delta$. Then
\[
    T_\perp^{p}\Omega\Delta=Z_\star\Sigma_\star^{p}\Omega_\star+Z_c\Sigma_c^{p}\Omega_c .
\]
Because $\Delta$ is nonsingular and $Z_\star^\top\Omega$ has full row rank, $\Omega_\star$ has full row rank, so that $\Omega_\star\Omega_\star^\dagger=I$, where $\Omega_\star^\dagger$ denotes the Moore--Penrose pseudoinverse. Moreover, $\mu_{k_\star+1}>\mu_{k_\star+2}\geq 0$ implies that $\Sigma_\star$ is nonsingular. Let
\[
    L=\Sigma_c^{p}\Omega_c\Omega_\star^\dagger\Sigma_\star^{-p},
    \qquad\text{so that}\qquad
    T_\perp^{p}\Omega\Delta\,\Omega_\star^\dagger\Sigma_\star^{-p}=Z_\star+Z_cL .
\]
Since $\Delta$ is nonsingular, $\range(Z_\star+Z_cL)\subseteq\range(T_\perp^{p}\Omega\Delta)=\range(T_\perp^{p}\Omega)=\range(Q)$. The matrix $Z_\star+Z_cL$ has full column rank, with $(Z_\star+Z_cL)^\top(Z_\star+Z_cL)=I+L^\top L$, and the orthogonal projector onto its range is
\[
    \Pi_L=(Z_\star+Z_cL)(I+L^\top L)^{-1}(Z_\star+Z_cL)^\top .
\]
Because $\range(\Pi_L)\subseteq\range(Q)$, we have $I-\Pi=(I-\Pi)(I-\Pi_L)$, and therefore
\[
    \norm{(I-\Pi)Z_\star}_2\leq\norm{(I-\Pi_L)Z_\star}_2 .
\]
Using $Z_\star^\top(Z_\star+Z_cL)=I$, we obtain
\[
    Z_\star^\top(I-\Pi_L)Z_\star=I-(I+L^\top L)^{-1}=L^\top L(I+L^\top L)^{-1},
\]
so that
\begin{equation}
\label{eq:tangent}
    \norm{(I-\Pi_L)Z_\star}_2^2=\frac{\norm{L}_2^2}{1+\norm{L}_2^2}\leq\norm{L}_2^2 .
\end{equation}
Finally, $\norm{\Sigma_c^{p}}_2=\mu_{k_\star+2}^{p}$, $\norm{\Sigma_\star^{-p}}_2=\mu_{k_\star+1}^{-p}$, and $\norm{\Omega_\star^\dagger}_2=1/\sigma_{\min}(\Omega_\star)$ give
\[
    \norm{L}_2\leq\left(\frac{\mu_{k_\star+2}}{\mu_{k_\star+1}}\right)^{p}\frac{\norm{\Omega_c}_2}{\sigma_{\min}(\Omega_\star)},
\]
which is the asserted bound.
\end{proof}

By \eqref{eq:tangent}, the right-hand side of \cref{lem:structural} bounds not only the sine but also the tangent of the largest principal angle between $\range(Z_\star)$ and $\range(Q)$. The same therefore holds for the bounds of \cref{thm:ds_gaussian_general,thm:ds_gaussian}. The tangent is the quantity reported in \cref{tab:illustration}.

\subsection{Gaussian degree-stratified sketches}

We use the following standard bounds for the extreme singular values of Gaussian matrices; see \cite[Theorem~II.13]{davidson2001local} and \cite[Corollary~5.35]{vershynin2012introduction}.

\begin{lemma}
\label{lem:gaussian_sv}
Let $G\in\mathbb R^{a\times b}$ have independent standard Gaussian entries, and let $t\geq 0$. Then
\[
    \mathbb P\bigl\{\norm{G}_2>\sqrt{a}+\sqrt{b}+t\bigr\}\leq e^{-t^2/2},
\]
and, if $a\leq b$,
\[
    \mathbb P\bigl\{\sigma_{\min}(G)<\sqrt{b}-\sqrt{a}-t\bigr\}\leq e^{-t^2/2}.
\]
\end{lemma}

We first state a version of \cref{thm:ds_gaussian} with explicit constants, in which a stratum that receives few sketch dimensions is admitted and does not contribute to the bound. We write $(x)_+=\max\{x,0\}$.

\begin{theorem}
\label{thm:ds_gaussian_general}
Let $Q$ be computed by \cref{alg:degree-stratified-rp} with Gaussian local sketches. For $t>0$, define
\begin{equation}
\label{eq:eta-phi}
    \eta_j(t)=\left(\frac{\bigl(\sqrt{\ell_j}-\sqrt{k_j}-t\bigr)_+}{\sqrt{n_j}+\sqrt{\ell_j}+t}\right)^{2},
    \qquad
    \Phi_t=\lambda_{\min}\Biggl(\sum_{j=1}^{s}\eta_j(t)\,W_j\Biggr).
\end{equation}
If $\Phi_t>0$, then, with probability at least $1-2s\,e^{-t^2/2}$,
\[
    \dist(Z_\star,Q)
    \leq
    \left(\frac{\mu_{k_\star+2}}{\mu_{k_\star+1}}\right)^{2q+1}
    \frac{1}{\sqrt{\Phi_t}} .
\]
\end{theorem}

\begin{proof}
Let $G_j=\sqrt{\ell_j}\,\Omega_j\in\mathbb R^{n_j\times\ell_j}$, which has independent standard Gaussian entries. For every $j$ with $k_j\geq 1$, let the columns of $\widehat V_j\in\mathbb R^{n_j\times k_j}$ form an orthonormal basis of $\range(V_j)$; then $V_j=\widehat V_j\widehat V_j^\top V_j$, and $\widehat V_j^\top G_j\in\mathbb R^{k_j\times\ell_j}$ has independent standard Gaussian entries. Consider the events
\[
    \mathcal E_j^{(1)}=\bigl\{\norm{G_j}_2\leq\sqrt{n_j}+\sqrt{\ell_j}+t\bigr\},
    \qquad j=1,\ldots,s,
\]
and, for every $j$ with $1\leq k_j\leq\ell_j$,
\[
    \mathcal E_j^{(2)}=\bigl\{\sigma_{\min}(\widehat V_j^\top G_j)\geq\sqrt{\ell_j}-\sqrt{k_j}-t\bigr\}.
\]
By \cref{lem:gaussian_sv}, each of these at most $2s$ events fails with probability at most $e^{-t^2/2}$, so that all of them hold simultaneously with probability at least $1-2s\,e^{-t^2/2}$. We assume in the following that they hold.

Let $\delta_j=(\sqrt{n_j}+\sqrt{\ell_j}+t)^{-1}$, and let $\Delta$ be the block-diagonal matrix with diagonal blocks $\delta_j\sqrt{\ell_j}\,I_{\ell_j}$, $j=1,\ldots,s$, so that
\[
    \Omega\Delta=[\delta_1P_1G_1,\ldots,\delta_sP_sG_s].
\]
We first bound $\norm{\Omega\Delta}_2$. Since $P_j^\top P_j=I$ and the strata are disjoint, we have $\sum_{j}\norm{P_j^\top\bs{x}}_2^2=\norm{\bs{x}}_2^2$ for every $\bs{x}\in\mathbb R^n$, and therefore
\[
    \norm{(\Omega\Delta)^\top\bs{x}}_2^2=\sum_{j=1}^{s}\delta_j^2\norm{G_j^\top P_j^\top\bs{x}}_2^2
    \leq\sum_{j=1}^{s}\delta_j^2\norm{G_j}_2^2\norm{P_j^\top\bs{x}}_2^2
    \leq\norm{\bs{x}}_2^2,
\]
where the last inequality uses the events $\mathcal E_j^{(1)}$. Hence $\norm{Z_c^\top\Omega\Delta}_2\leq\norm{\Omega\Delta}_2\leq 1$.

Next, we bound the smallest singular value of $Z_\star^\top\Omega\Delta=[\delta_1V_1^\top G_1,\ldots,\delta_sV_s^\top G_s]$ from below. For every $j$ with $1\leq k_j\leq\ell_j$, we use $V_j=\widehat V_j\widehat V_j^\top V_j$ and the event $\mathcal E_j^{(2)}$ to obtain
\begin{align*}
    V_j^\top G_jG_j^\top V_j
    &=(\widehat V_j^\top V_j)^\top(\widehat V_j^\top G_j)(\widehat V_j^\top G_j)^\top(\widehat V_j^\top V_j)\\
    &\succeq\sigma_{\min}(\widehat V_j^\top G_j)^2\,(\widehat V_j^\top V_j)^\top(\widehat V_j^\top V_j)
    \succeq\bigl(\sqrt{\ell_j}-\sqrt{k_j}-t\bigr)_+^2\,W_j ,
\end{align*}
where $(\widehat V_j^\top V_j)^\top(\widehat V_j^\top V_j)=V_j^\top V_j=W_j$. For the remaining indices, either $W_j=0$ or $(\sqrt{\ell_j}-\sqrt{k_j}-t)_+=0$, and the inequality $V_j^\top G_jG_j^\top V_j\succeq(\sqrt{\ell_j}-\sqrt{k_j}-t)_+^2\,W_j$ holds trivially. Summing over the strata, we obtain
\[
    (Z_\star^\top\Omega\Delta)(Z_\star^\top\Omega\Delta)^\top
    =\sum_{j=1}^{s}\delta_j^2\,V_j^\top G_jG_j^\top V_j
    \succeq\sum_{j=1}^{s}\eta_j(t)\,W_j ,
\]
and hence $\sigma_{\min}(Z_\star^\top\Omega\Delta)^2\geq\Phi_t>0$. In particular, $Z_\star^\top\Omega\Delta$, and therefore $Z_\star^\top\Omega$, has full row rank. The assertion now follows from \cref{lem:structural}.
\end{proof}

\begin{proof}[Proof of \cref{thm:ds_gaussian}]
Under the assumptions of the theorem, $t\leq\frac12\sqrt{\ell_j}$ and $k_j\leq k_\star$ give
\[
    \sqrt{\ell_j}-\sqrt{k_j}-t\geq\sqrt{\ell_j}-\sqrt{k_\star}-t\geq\tfrac12\sqrt{\ell_j},
    \qquad
    \sqrt{n_j}+\sqrt{\ell_j}+t\leq\sqrt{n_j}+\tfrac32\sqrt{\ell_j}\leq\tfrac52\sqrt{n_j},
\]
so that $\eta_j(t)\geq\ell_j/(25\,n_j)$ for all $j$. Since the matrices $W_j$ are positive semidefinite, it follows that
\[
    \sum_{j=1}^{s}\eta_j(t)\,W_j\succeq\frac{1}{25}\sum_{j=1}^{s}\frac{\ell_j}{n_j}\,W_j ,
\]
and hence $\Phi_t\geq\Phi/25$. Moreover, $\Phi\geq\min_j(\ell_j/n_j)>0$ because $\sum_jW_j=I$. The assertion follows from \cref{thm:ds_gaussian_general}.
\end{proof}

For a single stratum, $W_1=I$ and $k_1=k_\star$, and \cref{thm:ds_gaussian_general} reduces to the familiar estimate for the Gaussian range finder \cite{halko2011finding}, with $\Phi_t^{-1/2}=(\sqrt n+\sqrt\ell+t)/(\sqrt\ell-\sqrt{k_\star}-t)$. In general, $\eta_j(t)\approx\ell_j/n_j$ when $k_j\ll\ell_j\ll n_j$, so that $\Phi_t\approx\Phi$. The block structure enters the proof through the norm of the rescaled test matrix, which is the maximum of the block norms, each scaling like $\sqrt{n_j}$ rather than $\sqrt n$.

\subsection{The second-moment criterion}

The following proposition contains the statements on the second-moment criterion made in \cref{ssec:structural}.

\begin{proposition}
\label{prop:gram_mse}
Let $\Omega_1,\ldots,\Omega_s$ be Gaussian local sketches as in \cref{ssec:ds-sketch}. Let $Y=T_\perp\Omega$ with $\Omega$ as in \eqref{eq:structured-test-matrix}, and let $T_j=T_\perp P_j$. Then $\mathbb E[YY^\top]=T_\perp^2$ and
\begin{equation}
\label{eq:gram_mse}
    \mathbb E\left[\norm{YY^\top-T_\perp^2}_F^2\right]
    =\sum_{j=1}^{s}\frac{a_j^2+b_j^2}{\ell_j},
    \qquad
    a_j=\norm{T_j}_F^2,
    \quad
    b_j=\norm{T_jT_j^\top}_F .
\end{equation}
Moreover, for all positive numbers $\ell_1,\ldots,\ell_s$ with $\sum_{j}\ell_j=\ell$,
\begin{equation}
\label{eq:gram_mse_comparison}
    \sum_{j=1}^{s}\frac{a_j^2+b_j^2}{\ell_j}
    \geq\frac{1}{\ell}\Biggl(\sum_{j=1}^{s}\sqrt{a_j^2+b_j^2}\Biggr)^{2}
    \geq\frac{1}{\ell}\left(\norm{T_\perp}_F^4+\norm{T_\perp^2}_F^2\right),
\end{equation}
and the right-hand side is the value of \eqref{eq:gram_mse} for the standard Gaussian sketch, $s=1$.
\end{proposition}

\begin{proof}
Let $Y_j=T_j\Omega_j$, so that $YY^\top=\sum_jY_jY_j^\top$. Since $\mathbb E[\Omega_j\Omega_j^\top]=I_{n_j}$, we have $\mathbb E[Y_jY_j^\top]=T_jT_j^\top$, and since the strata form a partition of the node set, $\sum_jP_jP_j^\top=I$ and
\[
    \mathbb E[YY^\top]=\sum_{j=1}^{s}T_jT_j^\top=T_\perp\Biggl(\sum_{j=1}^{s}P_jP_j^\top\Biggr)T_\perp=T_\perp^2 .
\]
The matrices $Y_jY_j^\top-T_jT_j^\top$, $j=1,\ldots,s$, are independent and have mean zero, so that their Frobenius inner products have mean zero and
\[
    \mathbb E\left[\norm{YY^\top-T_\perp^2}_F^2\right]=\sum_{j=1}^{s}\mathbb E\left[\norm{Y_jY_j^\top-T_jT_j^\top}_F^2\right].
\]
For fixed $j$, write $\Omega_j=\ell_j^{-1/2}[\bs{g}_1,\ldots,\bs{g}_{\ell_j}]$ with independent standard Gaussian vectors in $\mathbb R^{n_j}$. Then $Y_jY_j^\top$ is the average of the $\ell_j$ independent and identically distributed matrices $\bs{x}_i\bs{x}_i^\top$ with $\bs{x}_i=T_j\bs{g}_i\sim\mathcal N(0,T_jT_j^\top)$, each with mean $T_jT_j^\top$, and therefore
\[
    \mathbb E\left[\norm{Y_jY_j^\top-T_jT_j^\top}_F^2\right]=\frac{1}{\ell_j}\,\mathbb E\left[\norm{\bs{x}\bs{x}^\top-T_jT_j^\top}_F^2\right],
    \qquad\bs{x}\sim\mathcal N(0,T_jT_j^\top).
\]
Since $\norm{\bs{x}\bs{x}^\top}_F^2=\norm{\bs{x}}_2^4$, the fourth-moment identity $\mathbb E\norm{\bs{x}}_2^4=\trace(T_jT_j^\top)^2+2\norm{T_jT_j^\top}_F^2$ yields
\[
    \mathbb E\left[\norm{\bs{x}\bs{x}^\top-T_jT_j^\top}_F^2\right]
    =\mathbb E\norm{\bs{x}}_2^4-\norm{T_jT_j^\top}_F^2
    =a_j^2+b_j^2 ,
\]
which proves \eqref{eq:gram_mse}.

The first inequality in \eqref{eq:gram_mse_comparison} is the Cauchy--Schwarz inequality,
\[
    \Biggl(\sum_{j=1}^{s}\sqrt{a_j^2+b_j^2}\Biggr)^{2}
    =\Biggl(\sum_{j=1}^{s}\sqrt{\frac{a_j^2+b_j^2}{\ell_j}}\,\sqrt{\ell_j}\Biggr)^{2}
    \leq\ell\sum_{j=1}^{s}\frac{a_j^2+b_j^2}{\ell_j}.
\]
For the second inequality, the Cauchy--Schwarz inequality in $\mathbb R^2$ gives $\sqrt{a_j^2+b_j^2}\sqrt{a_h^2+b_h^2}\geq a_ja_h+b_jb_h$ for all $j$ and $h$, and hence
\[
    \Biggl(\sum_{j=1}^{s}\sqrt{a_j^2+b_j^2}\Biggr)^{2}
    \geq\Biggl(\sum_{j=1}^{s}a_j\Biggr)^{2}+\Biggl(\sum_{j=1}^{s}b_j\Biggr)^{2}.
\]
Finally, $\sum_ja_j=\sum_j\norm{T_\perp P_j}_F^2=\norm{T_\perp}_F^2$, and the triangle inequality gives
\[
    \sum_{j=1}^{s}b_j=\sum_{j=1}^{s}\norm{T_jT_j^\top}_F\geq\Bigl\|\sum_{j=1}^{s}T_jT_j^\top\Bigr\|_F=\norm{T_\perp^2}_F .
\]
For $s=1$ we have $T_1=T_\perp$ up to a permutation of the columns and $\ell_1=\ell$, so that \eqref{eq:gram_mse} equals $(\norm{T_\perp}_F^4+\norm{T_\perp^2}_F^2)/\ell$.
\end{proof}

\subsection{Leverage and degree}

\begin{proof}[Proof of \cref{prop:dcsbm_leverage}]
Let $\beta_a=\sum_{b=1}^{k}B_{ab}\sum_{v_{i'}\in\mathcal C_b}\psi_{i'}$ for $a=1,\ldots,k$. For $v_i\in\mathcal C_a$, summing \eqref{eq:dcsbm} over $i'$ gives $d_i=\psi_i\beta_a$, and $\beta_a>0$ because the degrees are positive. Hence, for $v_i\in\mathcal C_a$ and $v_{i'}\in\mathcal C_b$,
\[
    S_{ii'}=\frac{\psi_i\psi_{i'}B_{ab}}{\sqrt{d_id_{i'}}}
    =\sqrt{\psi_i/\beta_a}\;B_{ab}\,\sqrt{\psi_{i'}/\beta_b} .
\]
Let $\nu_a=\sum_{v_i\in\mathcal C_a}\psi_i/\beta_a>0$, and let $J\in\mathbb R^{n\times k}$ have the entries
\[
    J_{ia}=\sqrt{\psi_i/(\beta_a\nu_a)}\quad\text{if } v_i\in\mathcal C_a,
    \qquad
    J_{ia}=0\quad\text{otherwise}.
\]
The columns of $J$ have disjoint supports and unit norms, so that $J$ has orthonormal columns, and
\[
    S=J\widetilde BJ^\top,
    \qquad
    \widetilde B_{ab}=\sqrt{\nu_a}\,B_{ab}\sqrt{\nu_b} .
\]
Since $B$ is nonsingular and $\nu_a>0$ for all $a$, the matrix $\widetilde B$ is nonsingular, so that $S$ has rank $k$ and $\range(S)=\range(J)$. The orthogonal projector onto $\range(S)$ is $JJ^\top$, and for $v_i\in\mathcal C_a$ its $i$th diagonal entry is
\[
    J_{ia}^2=\frac{\psi_i/\beta_a}{\nu_a}
    =\frac{\psi_i\beta_a}{\sum_{v_{i'}\in\mathcal C_a}\psi_{i'}\beta_a}
    =\frac{d_i}{\vol(\mathcal C(v_i))}.
\]

The trivial eigenvector $\bs{u}_1=D^{1/2}\bs{1}/\norm{D^{1/2}\bs{1}}_2$ satisfies $S\bs{u}_1=\bs{u}_1$ and therefore belongs to $\range(S)$; its entries satisfy $(\bs{u}_1)_i^2=d_i/\vol(\mathcal G)$. The eigenvectors of $S$ associated with the remaining nonzero eigenvalues span the orthogonal complement of $\bs{u}_1$ in $\range(S)$, which has dimension $k-1$ and orthogonal projector $JJ^\top-\bs{u}_1\bs{u}_1^\top$. The diagonal entries of this projector are the asserted leverage scores. If all communities have the same volume, then $\vol(\mathcal C_a)=\vol(\mathcal G)/k$ for all $a$, and hence $\zeta_i=(k-1)d_i/\vol(\mathcal G)$.
\end{proof}

The approximation \eqref{eq:gain} is obtained as follows. If $W_j=(\tau_j/k_\star)I$ and $\tau_j=k_\star m_j/\vol(\mathcal G)$ for all $j$, then
\[
    \Phi=\lambda_{\min}\Biggl(\sum_{j=1}^{s}\frac{\ell_j}{n_j}\,\frac{m_j}{\vol(\mathcal G)}\,I\Biggr)
    =\frac{1}{\vol(\mathcal G)}\sum_{j=1}^{s}\ell_j\,\frac{m_j}{n_j} .
\]
For the proportional allocation $\ell_j=\ell n_j/n$, the right-hand side equals $\ell\sum_jm_j/(n\vol(\mathcal G))=\ell/n$.

\subsection{Ritz vectors and spectral embedding}

\begin{proof}[Proof of \cref{thm:rayleigh_ritz}]
If $k=\ell$, then $U_kU_k^\top=\Pi$ and there is nothing to prove. Let $k<\ell$, and partition $\widetilde U=[\widetilde U_1,\widetilde U_2]$ with $\widetilde U_1\in\mathbb R^{\ell\times k}$, so that $U_k=Q\widetilde U_1$. Let $U_c=Q\widetilde U_2$ and $\Theta_c=\diag(\theta_{k+1},\ldots,\theta_\ell)$. Then $\Pi=U_kU_k^\top+U_cU_c^\top$ and
\[
    (I-U_kU_k^\top)Z_\star=(I-\Pi)Z_\star+U_cU_c^\top Z_\star .
\]
The two terms on the right-hand side have mutually orthogonal ranges, and hence
\begin{equation}
\label{eq:rr_split}
    \norm{(I-U_kU_k^\top)Z_\star}_2^2\leq\norm{(I-\Pi)Z_\star}_2^2+\norm{X}_2^2,
    \qquad X:=U_c^\top Z_\star .
\end{equation}

We derive a Sylvester equation for $X$. Since $C\widetilde U_2=\widetilde U_2\Theta_c$,
\[
    U_c^\top T_\perp\Pi=\widetilde U_2^\top Q^\top T_\perp QQ^\top=\widetilde U_2^\top CQ^\top=\Theta_cU_c^\top .
\]
Using $T_\perp Z_\star=Z_\star\Sigma_\star$ and $I=\Pi+(I-\Pi)$, we obtain
\[
    X\Sigma_\star=U_c^\top T_\perp Z_\star=\Theta_cX+E,
    \qquad
    E:=U_c^\top T_\perp(I-\Pi)Z_\star,
\]
that is, $X\Sigma_\star-\Theta_cX=E$. The eigenvalues of $\Sigma_\star$ are at least $\mu_{k_\star+1}$. By the Cauchy interlacing theorem \cite{parlett1998symmetric}, $\theta_{k+1}$ is at most the $(k+1)$st largest eigenvalue of $T_\perp$, which is at most $\mu_{k_\star+2}$ because $k\geq k_\star$. Hence the eigenvalues of $\Theta_c$ are at most $\mu_{k_\star+2}=\mu_{k_\star+1}-\gamma_\star$. Subtracting $\mu_{k_\star+1}X$ on both sides of the equation, we therefore have
\[
    X(\Sigma_\star-\mu_{k_\star+1}I)+(\mu_{k_\star+1}I-\Theta_c)X=E,
    \qquad
    \Sigma_\star-\mu_{k_\star+1}I\succeq 0,
    \qquad
    \mu_{k_\star+1}I-\Theta_c\succeq\gamma_\star I .
\]
Since $\norm{e^{-\xi(\mu_{k_\star+1}I-\Theta_c)}}_2\leq e^{-\xi\gamma_\star}$ and $\norm{e^{-\xi(\Sigma_\star-\mu_{k_\star+1}I)}}_2\leq 1$ for $\xi\geq 0$, the solution of this equation is
\[
    X=\int_0^\infty e^{-\xi(\mu_{k_\star+1}I-\Theta_c)}\,E\,e^{-\xi(\Sigma_\star-\mu_{k_\star+1}I)}\,d\xi,
\]
so that $\norm{X}_2\leq\norm{E}_2\int_0^\infty e^{-\xi\gamma_\star}\,d\xi=\norm{E}_2/\gamma_\star$; see also \cite{davis1970rotation} and \cite[Section~VII.2]{bhatia1997matrix}.

It remains to bound $E$. Since $U_c^\top=U_c^\top\Pi$ and $(I-\Pi)^2=I-\Pi$,
\[
    \norm{E}_2\leq\norm{\Pi T_\perp(I-\Pi)}_2\,\norm{(I-\Pi)Z_\star}_2 .
\]
The eigenvalues of $T_\perp$ lie in $[0,\mu_2]$, so that $\norm{T_\perp-\frac{\mu_2}{2}I}_2\leq\frac{\mu_2}{2}$, and because $\Pi(I-\Pi)=0$,
\[
    \norm{\Pi T_\perp(I-\Pi)}_2=\norm{\Pi\bigl(T_\perp-\tfrac{\mu_2}{2}I\bigr)(I-\Pi)}_2\leq\frac{\mu_2}{2}.
\]
Therefore $\norm{X}_2\leq\frac{\mu_2}{2\gamma_\star}\norm{(I-\Pi)Z_\star}_2$, and inserting this bound into \eqref{eq:rr_split} completes the proof.
\end{proof}

\begin{proof}[Proof of \cref{cor:main}]
The first statement follows by combining \cref{thm:ds_gaussian,thm:rayleigh_ritz}. For the second statement, let $\sigma_1,\ldots,\sigma_{k_\star}$ denote the singular values of $U_k^\top Z_\star\in\mathbb R^{k\times k_\star}$. They lie in $[0,1]$ and are the cosines of the principal angles between $\range(Z_\star)$ and $\range(U_k)$, so that
\[
    \norm{(I-U_kU_k^\top)Z_\star}_2^2=\lambda_{\max}\bigl(I-Z_\star^\top U_kU_k^\top Z_\star\bigr)=\max_i\,(1-\sigma_i^2).
\]
By the polar decomposition, there is a matrix $R\in\mathbb R^{k\times k_\star}$ with orthonormal columns for which $R^\top U_k^\top Z_\star$ is symmetric positive semidefinite with eigenvalues $\sigma_1,\ldots,\sigma_{k_\star}$. Hence
\begin{align*}
    \norm{U_kR-Z_\star}_F^2
    &=2k_\star-2\trace\bigl(R^\top U_k^\top Z_\star\bigr)
    =2\sum_{i=1}^{k_\star}(1-\sigma_i)\\
    &\leq 2\sum_{i=1}^{k_\star}(1-\sigma_i^2)
    \leq 2k_\star\,\norm{(I-U_kU_k^\top)Z_\star}_2^2 ,
\end{align*}
where we used $1-\sigma_i\leq 1-\sigma_i^2$ for $\sigma_i\in[0,1]$.
\end{proof}

\section*{Acknowledgments}
The authors used ChatGPT (OpenAI) for language editing
and to assist in reviewing mathematical arguments.
The authors assume responsibility for all content.

\section*{Data and code availability}
All datasets are publicly available from the cited sources, and the code is available from the corresponding author upon reasonable request.

\bibliographystyle{siamplain}
\bibliography{ref}
\end{document}

%% file: shared.tex
\usepackage{amsfonts,amssymb}
\usepackage{graphicx}
\usepackage{multirow}
\usepackage{booktabs}
\usepackage{siunitx}
\usepackage{algpseudocode}
\usepackage{enumitem}
\usepackage{placeins}

\DeclareGraphicsExtensions{.pdf,.png,.jpg}

\newsiamremark{remark}{Remark}
\crefname{remark}{Remark}{Remarks}
\AddToHook{env/remark/begin}{\crefalias{theorem}{remark}}

\algrenewcommand\algorithmicrequire{\textbf{Input:}}
\algrenewcommand\algorithmicensure{\textbf{Output:}}

\newcommand{\bs}[1]{\boldsymbol{#1}}
\newcommand{\norm}[1]{\left\lVert#1\right\rVert}
\newcommand{\orth}{\operatorname{orth}}
\newcommand{\range}{\operatorname{range}}
\newcommand{\diag}{\operatorname{diag}}
\newcommand{\rank}{\operatorname{rank}}
\newcommand{\trace}{\operatorname{trace}}
\newcommand{\vol}{\operatorname{vol}}
\newcommand{\dist}{\operatorname{dist}}
\newcommand{\nnz}{\operatorname{nnz}}

\headers{Degree-stratified randomized projection}{J. Eom, Y. Hwang, E. Kim, and Y. Jang}

\title{Degree-stratified randomized projection for spectral clustering%
\thanks{{\color{header1}\bfseries Funding:} This work was supported by the Global-Learning \& Academic Research
Institution for Master's and PhD Students, and Postdocs (LAMP) Program
(No.~RS-2024-00442775) and the BK21 FOUR (Fostering Outstanding Universities
for Research) Program (No.~5120200913674), funded by the Ministry of Education
and the National Research Foundation of Korea (NRF). This work was also
supported by the Regional Innovation System \& Education (RISE) Glocal
University 30 Program (No.~2026-RISE(Glocal University 30)-05-011), funded by
the Ministry of Education and Gwangju Metropolitan City through the Gwangju
RISE Center.}}

\author{
Junghyun Eom\thanks{Department of Mathematics and Statistics, Chonnam National
University, Gwangju, South Korea (\email{djawjdgus302@chonnam.ac.kr},
\email{youngjin@chonnam.ac.kr}, \email{flairz@chonnam.ac.kr},
\email{yongseok.jang@chonnam.ac.kr}). J.~Eom and Y.~Hwang contributed equally
to this work. Corresponding author: Y.~Jang.}
\and
Youngjin Hwang\footnotemark[2]
\and
Euihyun Kim\footnotemark[2]
\and
Yongseok Jang\footnotemark[2]
}